\documentclass [10pt]{article}
\usepackage{graphicx} 
\newtheorem{corollary}{Corollary}
\newtheorem{theorem}[corollary]{Theorem}

\newtheorem{conjecture}{Conjecture}
\newtheorem{proposition}[corollary]{Proposition}
\newtheorem{question}{Question}

\begin{document}

\title{Conformal restriction and related questions}
\date {Lecture Notes,  ICMS Edinburgh, July 2003}
\author{Wendelin Werner}
\maketitle
\centerline{Universit\'e Paris-Sud and IUF}

\centerline {Laboratoire de Math\'ematiques, Universit\'e Paris-Sud}

\centerline{B\^at. 425, 91405 Orsay cedex, France}

\centerline {\textit{e-mail: wendelin.werner@math.u-psud.fr}}

\centerline {http://www.math.u-psud.fr/$\sim$werner}

\hfuzz =3pt
\
\newcommand{\Prob} {{\bf P}}
\font \m=msbm10
\newcommand{\R}{{\hbox {\m R}}}
\newcommand{\C}{{\hbox {\m C}}}
\newcommand{\Z}{{\hbox {\m Z}}}
\newcommand{\N}{{\hbox {\m N}}}
\newcommand{\U}{{\hbox {\m U}}}
\def\H{{\hbox {\m H}}}
\def\P{{\bf P}}
\def \expect {{\bf E}}
\def \eps {\varepsilon}

\vskip 20mm

These notes are based on mini-courses given in July 2003 at the University of St-Andrews, and at the ICMS in 
Edinburgh. The goal of these lectures is to give a self-contained sketchy and heuristic survey of the
recent results concerning conformal restriction, that were initiated in our 
joint work with Greg Lawler and Oded Schramm \cite {LSWr,LSWSAW}, and further investigated in the 
papers \cite {FW, Whid, LWbls, Dub}. It also finds some of its roots in the earlier ``pre-SLE''
paper \cite {LW2}.
If the reader wants full proofs and cleaner statements, (s)he is advised to consult these original papers. We will also not go into the details of the construction of the Schramm-Loewner Evolution (SLE), referring the motivated and interested reader to \cite {Wstf,Lbook} and the references therein for more details. These notes can be viewed as complementary to my Saint-Flour  notes \cite {Wstf} (one can read either of them before the other), and are maybe aimed at a less probabilistic audience.  

The papers dealing with the restriction property all take the existence and properties of SLE for granted, and are therefore maybe difficult to read if one is not already acquainted to the SLE background. We will try in these lectures to conversely use the intuition and ideas building on the restriction properties and self-avoiding walks to shed some light onto some of the properties of SLE, and its relation to some ideas from conformal field theory. We hope that this can be useful in the theoretical physics community as well. 

Almost all results in these notes are borrowed from the existing above-mentioned papers, but there are a few new ones, such as for instance the construction of the one-sided 
restriction measures via Poissonian clouds of Brownian excursions.

\vskip 2cm

{\bf Acknowledgements.}
Thanks are of course due to Greg Lawler and Oded Schramm. I would also like to thank Vincent Beffara, Roland Friedrich and Julien Dub\'edat, whose PhD's I have (or had) the privilege of supervising for the many fruitful (still ongoing) discussions. 
\newpage

\tableofcontents

\newpage
 
\section {Introduction}

\subsection {Generalities}

One of the main issues in probability theory and statistical physics is to understand the large scale behaviour of 
random systems, that are often defined in the discrete setting. For instance, one has a well-defined probability measure
on a finite state-space, and one lets the size of the state-space go to infinity, and tries to understand the 
asymptotic behaviour of some observables. In many cases, the asymptotic behaviour is deterministic, but it can 
also happen to be random. The existence of the scaling limit is usually justified heuristically via a renormalization 
(or fixed point) argument. But, in the generic case,
 a further mathematical description is out of reach. The complexity of the 
system can make it impossible to encode the randomness in a proper mathematical way.
In some exceptional cases, an additional (combinatorial, algebraic, analytical) feature 
can be shown to hold. This extra structure can then be used to show pin down this scaling limit 
and to encode this complexity. It can give rise to unusual random processes (where for instance the noise 
in the sense of Tsirelson \cite {Tsi} is not necessarily white), but that can be (precisely because of this 
additional mathematical structure) related and of interest to other areas of mathematics.

Two-dimensional critical systems are believed to belong to this class. The additional structure (that has been proved to hold in some cases \cite {Sm,LSWSAW}) is conformal invariance. This had been recognized long ago by theoretical 
physicists (see e.g. \cite {BPZ}) and gave rise to an intense activity for instance in conformal field theory.
More recently, the SLE (SLE stands for Schramm-Loewner Evolutions) approach did provide a simple mathematical new perspective to these systems. In the present lectures, we will focus on one rather specific aspect of these random systems that we initially thought of because of the problems of self-avoiding walks, but turned out in the end to be relevant  to all these systems. The global idea, which  recalls some considerations of conformal field theory, is basically
to see how the law of these random curves behave when one changes (or perturbs infinitesimally) the domain it is defined in. 

\subsection {Brownian motion, conformal invariance}

Suppose that we are looking for a ``uniform'' probability measure on the space of $d$-dimensional 
continuous curves. The state-space is infinite, so that this notion is rather vague,
but it is easy to see that the natural candidate for such a measure is 
$d$-dimensional Brownian motion. One standard way to proceed is to start with a discretization of the 
state-space: For any finite $N$, the uniform measure 
on the set of paths of length $N$ on a given lattice and fixed starting point is just 
the law of simple random walk on that lattice. The continuous limit of simple random 
walk is (under mild conditions) Brownian motion, regardless of the precise lattice that one starts with.

It is also possible to restrict the class of paths. For instance, one can consider a 
finite domain $D$, a point $O$ in its interior, and try to construct a measure on the 
set of paths from $O$ to the boundary of the domain. Let us consider the measure on 
Brownian paths, that are started from $O$ and stopped at their first exist of $D$.
Note that, even though this measure is closely related to the previous ``uniform measure'', it is 
not ``uniformly distributed'' among the family of paths from the origin to $\partial D$:
In the discrete case (on the square lattice in two dimensions, say), if one considers the law of random walk started from $O$ and stopped at its first hitting of the boundary of the domain, it assigns a 
probability of $4^{-n}$ to each admissible path of $n$ steps, and $n$ is varying from one path to another.
One ``penalizes'' the mass of a path according to its length.

As pointed out by Paul L\'evy \cite {Le}, planar Brownian motion is conformally invariant. 
This means that if one considers a planar Brownian $Z$ started from $O$ and stopped at its 
first exit time $T$ of a simply connected domain $D$, and if $\Phi$ denotes a conformal map from $D$ onto some other domain $\tilde D$ (i.e. a one-to-one map that preserves angles), then the law of $\Phi (Z)$ is that of a Brownian motion 
started from $\Phi(O)$ and stopped at its first exit of $\tilde D$. Actually, this is not completely true, because one 
has to reparametrize time in a proper way. In fact, the rigorous statement is that 
$$ 
\Phi (Z_t) = \tilde Z_{\int_0^t | \Phi' (Z_s)|^2 ds},
$$
for all $t \le T$, where $\tilde Z$ is a Brownian motion started from $\Phi (O)$, that is stopped at 
$\tilde T = \int_0^T |\Phi' (Z_s)|^2 ds$, which is its exit time of $\tilde D$.
For instance, if $\Phi (z)= 2z$, then one has to speed up time by a factor of 4, so that 
$2 Z_t$ is in fact Brownian motion running at speed $4t$.

This shows that in general, the image of a Brownian path of fixed prescribed time-length under a conformal 
transformation is not a Brownian path with fixed and prescribed time-length. The ``uniform'' distribution 
is not fully preserved under conformal transformation. In this respect, it is more natural to deal with 
the law of Brownian motion with given endpoints or stopped at stopping times.
For instance, if $P_{D,O}$ denotes the law of Brownian motion started from $O$ and stopped at its 
first exit of $D$, we see that modulo time-reparametrization, $\Phi \circ P_{D,O}$ (which means the image measure of 
$P_{D,O}$ under the mapping $Z \mapsto \Phi (Z)$) is identical to $P_{\Phi (D), \Phi (O)}$.

One could work with such paths from an inner point to the boundary of a domain (and this would give rise 
to the ``radial restriction'' theory \cite {LSWrr}), but we will in this paper only speak about paths from one 
boundary point to another boundary point of a domain.
It is possible to define the natural Brownian measure on paths from one point $A$ of the boundary of 
$D$ to another point $B$ on the boundary of the domain $D$. 
In the discrete case on the square lattice, the mass of a path is 
again proportional to $4^{-n}$,
where the renormalizing constant corresponds to the conditioning 
of the random walk by the event that it exits $D$ at $B$ when starting from $A$.
In the scaling limit, this process can be understood as Brownian motion started from $A$ 
and conditioned to exit $D$ at $B$. Even if this conditioning does a priori not make sense 
(since $A \in \partial D$, it is an event of zero probability), it is not difficult to 
make this rigorous (for instance, using an $h$-process, or by
letting the starting point tend to $A$ from the inside of the domain, and condition the Brownian motion 
to exit in a neighbourhood of $B$).
We will call $P_{D,A,B}=P^{BM}_{D,A,B}$ the law of this process (and drop the superscript when there is 
no ambiguity). Conformal invariance then also holds for these processes:

\begin {proposition}
\label {p.ci}
If $\Phi$ is a conformal transformation from $D$ onto another domain $\tilde D$, and
if the law of $Z$ is $P_{D,A,B}$, then the law of $\Phi (Z)$ is 
$P_{\Phi(D), \Phi(A), \Phi (B) }$, modulo increasing reparametrization 
of the path.
\end {proposition}

\subsection {Restriction}

Suppose that $D' \subset D$ are simply connected (and $D \not= \C$), 
and that $A,B \in \partial D \cap \partial D'$.
In the discrete case, consider the law $P_D$ (respectively $P_{D'}$) of a simple random walk
$\omega$  started
from $A$ and conditioned to exit $D$ at $B$ (resp. to exit $D'$ at $B$).
Clearly, if $\omega$ is sampled according to $P_D$, but conditioned on the event that $\{\omega \subset D'\}$,
then the resulting law is exactly $P_{D'}$.
This is basically due to the fact that the probability of the path $\omega$ is proportional to $4^{-n}$,
whether it lives in $D$ or $D'$. The ``energy'' of the path ($n \log 4$ here, so that the 
probability is proportional to $\exp (-\hbox {energy})$ ) is depending only on $\omega$ itself, and not on the space of paths one considers. 
This property still holds in the scaling limit:
For any $D' \subset D$, for any $A$ and $B$ on $\partial D \cap \partial D'$,
$$
P_{D,A,B}^{BM} \mid \{ Z \subset D'\} = P_{D',A,B}^{BM}.
$$
We call this the restriction property.

Note that if we know $P_{D,A,B}^{BM}$, we have now two ways to get $P_{D',A,B}^{BM}$ for free.
The first one is by conditioning, the second one by conformal invariance (choosing a conformal map
$\Phi$ from $D$ onto $D'$ that leaves $A$ and $B$ invariant -- such maps exist by Riemann's mapping theorem).
A priori, one might wonder whether it is at all possible to find a measure on paths $P_{D,A,B}$ such that these
two ways coincide (for all $D' \subset D$), but we have just seen that the Brownian measure does the job.
On the other hand, this condition seems quite strong, so that one can ask if there exist other
measures on paths that also satisfy it. This will be one of the main issues in these lectures.

This leads to the following abstract definition:
Assume that a family of measures $P_{D,A,B}$ on curves $\omega$ from $A$ to $B$ in $D$
(that is indexed by $(D,A,B)$) satisfies:

\begin {itemize}
\item For all open $D' \subset D$, and $A \not= B$ on $\partial D \cap \partial D'$: 
The measure $P_{D,A,B}$ conditioned on $\{\omega \subset D'\}$ is the 
probability measure $P_{D',A,B}$.
\item
For any conformal transformation $\Phi$ on $D$, the image of $P_{D,A,B} $ under 
$\Phi$ is $P_{\Phi(D), \Phi(A), \Phi(B)}$ (modulo time-change).
\end {itemize}

We then say that this family satisfies {\em conformal restriction}. Let is again insist on the 
fact that this is a rather strong condition: Conformal invariance basically shows that all $P_{D,A,B}$ are 
defined from just one of them (for instance $P_{\H,0,\infty}$). 
Restriction then gives an additional relation between all these measures. 

An alternative rephrasing of the conformal restriction in terms of one measure 
$P_{D,A,B}$ (here $P_{\H,0,\infty}$) goes as follows:
Suppose that the family $P_{D,A,B}$ satisfies conformal restriction, then let $\gamma$
be a random path with law $P_{\H, 0 , \infty}$. Then:
\begin {enumerate}
\item
For any $\lambda > 0$, the law of $\lambda \gamma$ is equal to the law of $\gamma$
(modulo time-reparametrization).
\item
For any $H \subset \H$ that has the origin and infinity on its boundary (we can in fact also assume for convenience that $\H \setminus H$ is bounded and bounded away from infinity),
define the conformal map $\Phi_H$ from $H$ onto $\H$ that preserves the origin and such that $\Phi_\H (z) \sim z$ when $z \to \infty$ (this map exists and is unique, by Riemann's mapping Theorem). Then, the
law of $\gamma$ conditioned on $\{ \gamma \subset H \}$ is identical to the law of 
$\Phi_H^{-1} (\gamma)$. 
In words, this means that the law of $\gamma$ conditioned to remain in $H$ is identical to 
the law of the conformal image $\Phi_H^{-1} (\H)$ of the curve.
\end {enumerate}
The first fact follows from the fact that $z \mapsto \lambda z$ is a conformal transformation from
$\H$ onto itself.

Conversely, if a random curve $\gamma$ from the origin to infinity in $\H$ satisfies these two conditions, then
one can define for all $(D,A,B)$, the law $P_{D,A,B}$ of $\Phi(\gamma)$, where $\Phi$ is a conformal transformation
from $\H$ onto $D$ with $\Phi (O)=A$ and $\Phi (\infty)= B$, and check that this family of laws satisfies conformal 
restriction.
Therefore if these two conditions hold, we will sometimes say that the law of $\gamma$ satisfies conformal restriction.

\subsection {Motivation from self-avoiding walks}

One of our initial motivations was to reach a better understanding of self-avoiding walks.
In the discrete case, this is the uniform measure on the set of paths of length $N$ with a given starting point on a given lattice. 

When $N$ grows to infinity, it is easy by sub-multiplicativity, to see that the number $a_N$
of such self-avoiding walks of length $N$ on a given lattice grows at first order exponentially with $N$. More precisely, since $a_{N+M} \le a_N a_M$ and $a_N \ge 2^N$ (see e.g. \cite {MS}), there exists a lattice-dependent positive constant $\mu$ such that $a_N^{1/N}$
converges to $\mu$ (recall that without the self-avoiding constraint, the number of walks of length $N$ on the square lattice is $4^N$) when $N \to \infty$. 
For the same reasons as before, in order to exploit conformal invariance,
it will be convenient
to fix the endpoints of the curve instead of its length. The natural attempt is therefore to 
consider the discrete measure on paths from $A$ to $B$ in $D$ (on a discrete lattice) 
that puts a weight proportional to $\mu^{-n}$ to a self-avoiding path 
from $A$ to $B$ in $D$ with $n$ steps, and to try to understand the limiting object, 
when the mesh-size of the lattice goes to zero.
Again, this gives an ``intrinsic'' measure on paths, such that if one conditions 
$P_{D'}$ to those curves that stay in $D$, one obtains $P_D$. But, as opposed to the 
random walk/Brownian case,  the 
existence of the limit when the mesh-size vanishes is still an open problem.

 While there are basically no rigorous mathematical 
 results concerning long self-avoiding curves, various 
 striking predictions have been formulated by theoretical physics. 
 Their arguments often invoke conformal field theory, and are at present not 
 well-understood on a rigorous mathematical level. We hope that the results surveyed in the 
 present notes gives to these predictions a clearer status.  
 Let us briefly mention one of their conjectures (see Nienhuis \cite {N}):
 
 When the mesh-size vanishes, the appropriately rescaled random self-avoiding curve should look like 
 continuous curves with fractal dimension $4/3$. This can be also formulated in the following weaker form: The typical diameter of a self-avoiding walk with $N$ steps is of the order of $N^{3/4}$ (in the same way as the typical diameter of a simple random walk with $N$ steps is roughly $N^{1/2}$ because of the central limit theorem).
This exponent $3/4$ had first been proposed by Flory, in the 1940's \cite {Fl}.
   
 Note that the notion of conformal invariance itself (that is implicitly used in conformal 
 field theory) has to our knowledge not be given a clean precise meaning in the theoretical physics 
 literature.
Here, in our setting with fixed endpoints, it can formulated as follows:

\begin {conjecture}
The scaling limits $P_{D,A,B}^{SAW}$ of the measures on self-avoiding curves from $A$ to $B$ in $D$
(that puts weight proportional to $\mu^{-n(\omega)}$ to each such walk $\omega$, where $n(\omega)$ is its number of steps) 
exist. Furthermore, there are conformally invariant in the sense that they are satisfy Proposition~\ref {p.ci}.
\end {conjecture}

Note that, just as the measures on simple random walks do, the discrete measures on self-avoiding curves 
satisfy restriction as well. This property should clearly be preserved in this scaling limit. Hence, if the
previous conjecture holds, then clearly,
the family $P_{D,A,B}^{SAW}$ should also satisfy conformal restriction.

This has to be compared with the following theorem from \cite {LSWr} (see \cite {LSW4/3} for the
dimension)  that we shall discuss in these lectures:

\begin {theorem}
\label {main1}
There exists a unique probability measure on continuous paths without double points, that satisfies conformal restriction. It is the chordal Schramm-Loewner Evolution (SLE) with parameter $8/3$, and it is supported on 
curves with fractal dimension $4/3$.
\end {theorem}

\begin{figure}
\centerline{\includegraphics*[height=3in]{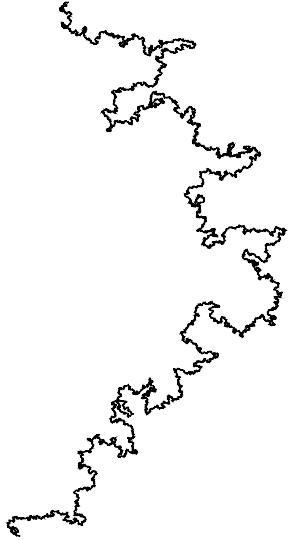}}
\caption{\label{f.saw}
Sample (courtesy of Tom Kennedy) of the beginning of an infinite half-plane walk (conjectured to 
converge to chordal $SLE_{8/3}$).}
\end{figure}

This of course leads to:

\begin {conjecture}
 The scaling limit of the measure on self-avoiding curves exists and it is this SLE$_{8/3}$.
\end {conjecture}

It is easy to construct other measures on curves that are ``intrinsic''. For instance, 
one can consider a measure on paths with double points but no triple points, or no self-crossings.
 Or measures that penalize paths according to their number of self-intersections (these are 
often called weakly self-avoiding walks in the literature). For each model (and each lattice), there 
exists a lattice dependent constant $\mu$ (that can be viewed as a critical value) such that the 
law of the limiting model should exist and exhibit interesting features.
If the limit exist and is conformally invariant (which is not always the case), then it should satisfy 
conformal restriction, and be related to the measures that we will be discussing in these lectures.

\subsection {Remarks on the half-plane measures}

\begin{figure}
\centerline{\includegraphics*[height=3in]{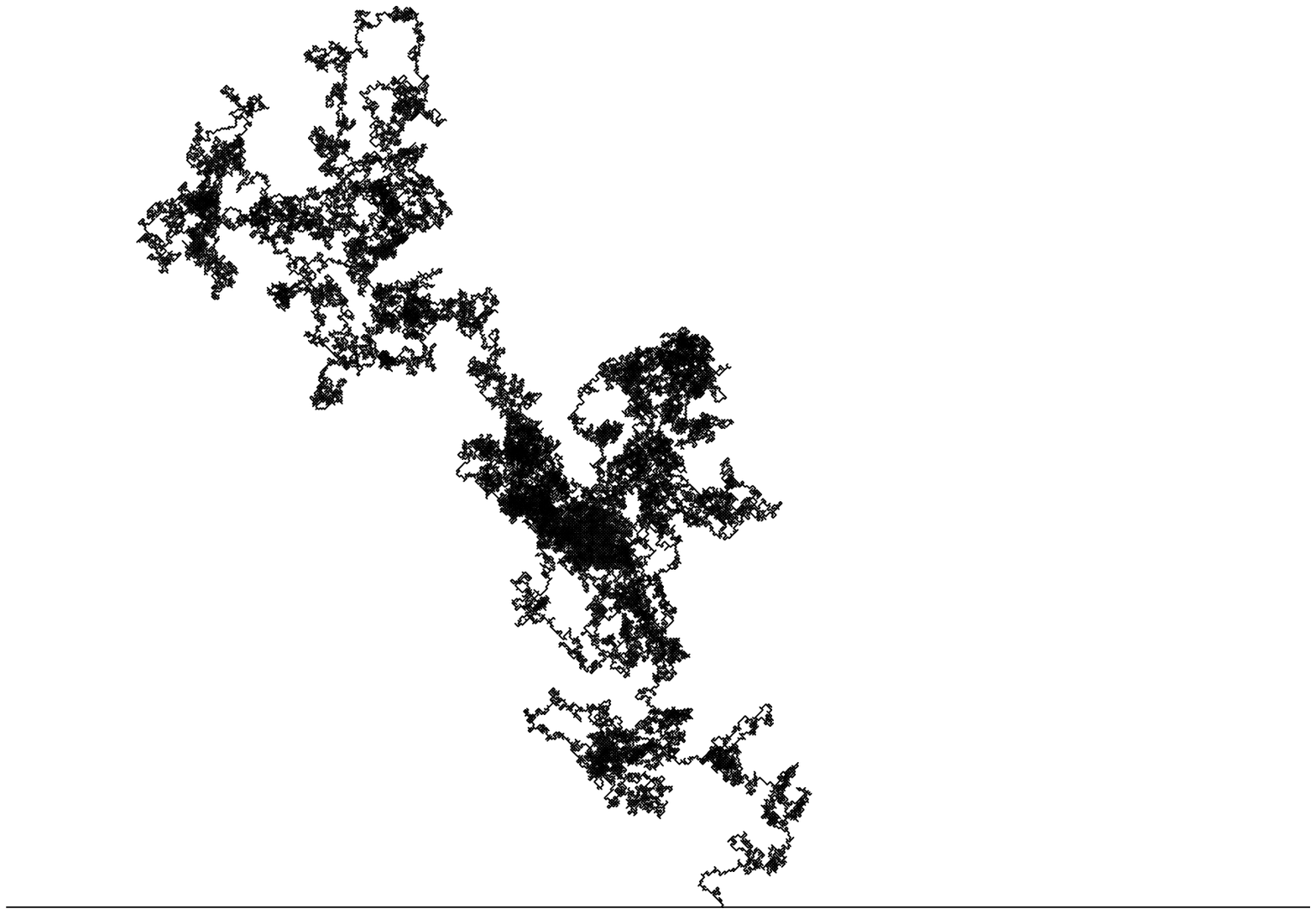}}
\caption{\label{f.ex}
Beginning of a sample of a discrete approximation of $P_{\H,0,\infty}^{BM}$.} 
\end{figure}

If one wishes to define the measure $P_{\H,0,\infty}^{BM}$, one can not use the 
same discrete approximation 
as before because the paths are infinite (so that $4^{-n} =0$ etc.). 
In the discrete square lattice (if one considers a continuous time Markov chain), one can
describe the corresponding random walk as follows:
 The real coordinate $X$ jumps like ordinary random walk (there is no conditioning),
while the imaginary part $Y$ is random walk ``conditioned to never hit zero''. In
other words, if $Y_t = y$, then at the next vertical jump, the walk moves up with probability $1 + 1/y$ and down with probability $1-1/y$.
The scaling limit of $Y$ is the three-dimensional Bessel process, which can be interpreted as Brownian motion conditioned to remain forever positive (and for this reason, 
this process tends almost surely to infinity when time grows to infinity), see e.g. \cite {RY}. 
Note that (both in the discrete and continuous picture), the law of 
the two-dimensional process $Z_t = (X_t, Y_t)$ in the upper half-plane can be understood as follows:
For any fixed $t$, the law of $(Z_s, s \le t)$ is the limit 
of that of the restriction to the time-interval $[0,t]$ 
of a planar Brownian motion (or simple random walk) of length $T$ conditioned to remain in the upper half-plane on the whole time-interval $[0,T]$, when $T \to \infty$. In this sense, 
$P_{\H, 0, \infty}^{BM}$ can be viewed as a uniform measure. This last remark also holds for the measures
$P_{H, 0, \infty}$, where $H \subset \H$ is a simply connected subset of $\H$ such that
$\H \setminus H$ is bounded and bounded away from the origin (in the sequel, when we use the notation $H$, we will 
implicitly mean such sets). 

It is also possible to define rigorously the measure on infinite discrete self-avoiding paths 
in the upper half-plane: One can consider the uniform measure on 
self-avoiding walks $S$ of length $N$ on $\Z^2$ 
started from the origin, and that stay in the 
upper-half plane. Then, when $n$ is fixed and $N$ tends to infinity, the law of 
$S(0), \ldots, S(n)$ can be proved to converge \cite {LSWSAW}, building on arguments of 
Kesten and  Madras-Slade \cite {Kes,MS}.
Since this holds for all $n$, this 
defines a law on infinite self-avoiding walks $S$ from $0$ to infinity in the upper half-plane, that can be 
heuristically understood as the uniform measure on infinite self-avoiding walks from the origin to infinity in the upper half-plane. In this case, the scaling limit problem is the existence of the limit of the law of 
the path of $\delta S$ when $\delta \to 0$ (corresponding to the walk on the lattice
$\delta \Z^2$). 

\section {The restriction exponent}
\label {S3}
\subsection {The exponent for the Brownian measure}

Let us first focus on the Brownian measure $P= P_{\H,0,\infty}^{BM}$.
The following description of $P$ will be useful:
Consider a Brownian path $Z$ that is started from $i \eps$, and that is 
conditioned to hit the line $\{\Im z = R \}$ before the real line. 
When $\eps \to 0$ and $R \to \infty$ (regardless of the order of the limits), 
the law of $Z$ converges to $P$ in some appropriate sense.
Note that the probability that a Brownian motion started from $i\eps$ hits
the line $\{\Im z = R \}$ before the real line is exactly $\eps/R$
(because $\Im Z$ is a martingale, i.e. its mean-value is constant). 
This corresponds to start $Z$ very very close to $0$ and to reach a ``neighbourhood'' of the boundary point 
at infinity.

Suppose that the law of $Z$ is $P$. Suppose that $\eps$ is very small, $R$ is very large. The law of $Z$ is 
close to $R/ \eps$ times the law of a Brownian motion started from
$i\eps$, restricted to the event that it hits $\{\Im (z) =R\}$ before 
the real line.
What is the probability that $Z$ also stays in $H$?
The previous description shows that it is 
 close to $R/ \eps$ times the probability that a Brownian motion 
started from $i\eps$ hits $\{\Im (z) =R\}$ before exiting $H$.

We now consider the image of this Brownian motion under 
the mapping $\Phi_H$ that is defined as follows:
$\Phi_H$ is the unique conformal map from $H$ onto $\H$ 
such that 
$$\Phi_H (0) =0 \hbox { and } 
\Phi_H (z) \sim z \hbox { when } z \to \infty$$
(by Riemann's mapping Theorem, this mapping exists and is unique).
We will use this definition of $\Phi_H$ throughout the paper.
If one looks at the image of this Brownian motion under the mapping $\Phi_H$,
one sees that this is $R / \eps$ times the probability that a Brownian motion 
started (near) from $i \Phi'(0) \eps$ hits $\Im (z) = R$ (since $\Phi(z) 
\sim z$ at infinity, this is close to $\Im (\Phi^{-1}(z)) = R$) before the 
real line. In other words, this probability is close to 
$\Phi_H'(0) \eps / R \times (R/ \eps )= \Phi_H' (0)$.
In the limit when $\eps \to 0$ and $R \to \infty$, we see that
\cite {V,LSWr} for all $H$,
\begin {equation}
\label {1}
P [ \gamma \subset H ] = \Phi_H' (0)
.\end {equation}
Here (and in the sequel), when we write $\gamma \subset H$, we will implicitly mean that
$\gamma \subset H \cup \{ 0 \}$.

\subsection {Restriction exponent}

As we shall now see, the fact that this probability in (\ref {1})
 is a power of $\Phi_H'(0)$ is 
in fact a general feature of conformal restriction. This power ($1$, in the Brownian case) will be 
called the restriction exponent.

Suppose that a random 
curve $\gamma$ from the origin to infinity in the upper half-plane satisfies conformal restriction. Note that it 
implies that its law is scale-invariant (because $z \mapsto \lambda z$ is a conformal transformation that
maps $\H$ onto itself). One can view the probability $P[\gamma \subset H] $ as a function 
$f( \Phi_H)$ of $\Phi_H$.
Recall that the law of $\Phi_H(\gamma)$, when $\gamma$ is conditioned to stay in $H$,
is identical to the law of $\gamma$ itself, so that the probability that it 
stays in some other set $H'$ is $f(\Phi_{H'})$. Hence,
\begin {eqnarray*}
f ( \Phi_{H'}\circ \Phi_H )
&=& P [ \gamma \subset \Phi_H^{-1} \circ \Phi_{H'}^{-1} ( \H ) ] \\
&=& P [ \gamma \subset \Phi_H^{-1} \circ \Phi_{H'}^{-1} ( \H )  \mid \gamma \subset H]
P [ \gamma \subset H ] \\
&=& f(\Phi_H) \times f(\Phi_{H'}). 
\end {eqnarray*}
In other words, the function $f$ is a homomorphism from the semi-group of conformal 
mappings $\Phi_H$ onto the multiplicative semi-group $[0,1]$. 
We will now briefly and heuristically
 justify that it implies in fact that $f( \Phi_H) = \Phi_H'(0)^\alpha$
for some exponent $\alpha$:

Loewner's theory (and we will come back to this later)
shows that there exists a one-parameter family $\varphi^t$ of 
mappings (in fact $\varphi^t= \Phi_{\H \setminus \eta [0,t]}$ for a well-defined 
curve $\eta$) such that 
$$\varphi^{t+s} = \varphi^t \circ \varphi^s$$
 for all $t, s$.
It follows immediately that for some constants $c$ and $\alpha$, 
$$(\varphi^t)'(0)  = \exp (-ct) \hbox { and } f(\varphi^t) = \exp ( -c \alpha t) = (\varphi^t)'(0)^\alpha.
$$

If we define for any positive real $x$, $\varphi_x^t (z) = x \varphi^t (z/x)$, we see 
immediately that $(\varphi_x^t)' (0) = (\varphi^t)'(0)$, and the scale-invariance of $\gamma$
 implies that 
$f(\varphi_x^t) = f( \varphi^t)$.
Similarly, one can see that $f$ (and the derivative at $0$) are invariant under conjugation 
with respect to the symmetry $\sigma$ with respect to the imaginary axis.

But Loewner's theory (e.g. \cite {Dur})  shows
that it is fact possible to approximate any mapping $\Phi_H$ 
by the iteration of many conformal maps $\phi_1, \ldots , \phi_n$ such that each 
$\phi_j$ is a conformal map of the type $\varphi_x^t$ (or $\sigma \circ \varphi_x^t \circ \sigma$).
It follows that 
\begin {eqnarray*}
f( \Phi_H)
&\sim&
f( \phi_1 \circ \cdots \circ\phi_n) \\
&=&
f(\phi_1) \times \cdots \times f(\phi_n) \\
&=&
\phi_1'(0)^\alpha \times \cdots \times \phi_n'(0)^\alpha \\
&=& (\phi_1 \circ \cdots \circ \phi_n)'(0)^\alpha \\
& \sim&
\Phi_H'(0)^\alpha.
\end {eqnarray*}
We have just indicated a rough justification of the fact that:

\begin {proposition}
If the path $\gamma$ satisfies conformal restriction, then there exists a constant 
$\alpha>0$ such that for all $H$,
\begin {equation}
\label {eq2}
P [\gamma \subset H ] = \Phi_H'(0)^\alpha.
\end {equation}
We call $\alpha$ the restriction exponent of $\gamma$.
\end {proposition}

\subsection {Characterization of the filling  by the exponent}

Suppose for a moment that $\gamma$ is a random simple curve (i.e. with no double points)
from the origin to infinity in the upper half-plane that satisfies conformal restriction.
We have just seen that there exists a positive constant $\alpha$ such that for all
$H$, (\ref {eq2}) holds. Furthermore, since $\gamma$ is a simple curve, it is in fact 
not difficult to see that the knowledge of $P[\gamma \subset H]$ for all $H$ characterizes the law of 
$\gamma$ (modulo time-reparametrization). Hence, this reduces the possible laws of such random curves $\gamma$
to a one-dimensional family indexed by the parameter $\alpha$. We shall later see that 
in fact, only one value of $\alpha$ gives rise to the law of a random simple curve. 

Let us now see what the formula (\ref {eq2}) tells about the law of 
$\gamma$, if we do not assume a priori that $\gamma$ is a simple curve.
Define the filling ${\cal F}(\gamma)$ of a curve $\gamma$ as the set of points 
in $\H \setminus \gamma$ that are not in the connected components of $\H \setminus 
\gamma$ that have $(0, \infty)$ and $(- \infty, 0)$ on their boundaries.
If $\gamma$ is a continuous curve from $0$ to infinity in $\H$,
then ${\cal F} (\gamma)$  is a closed connected subset of $\H$, and 
$\H \setminus {\cal F} (\gamma)$ consists of two unbounded connected 
components: $C_+ (\gamma)$ and $C_- (\gamma)$ that have respectively $\R_+$ and
$\R_-$ as parts of their boundaries.
When $\gamma$ is a simple curve, then ${\cal F} (\gamma) = \gamma$.

\begin{figure}
\centerline{\includegraphics*[height=2in]{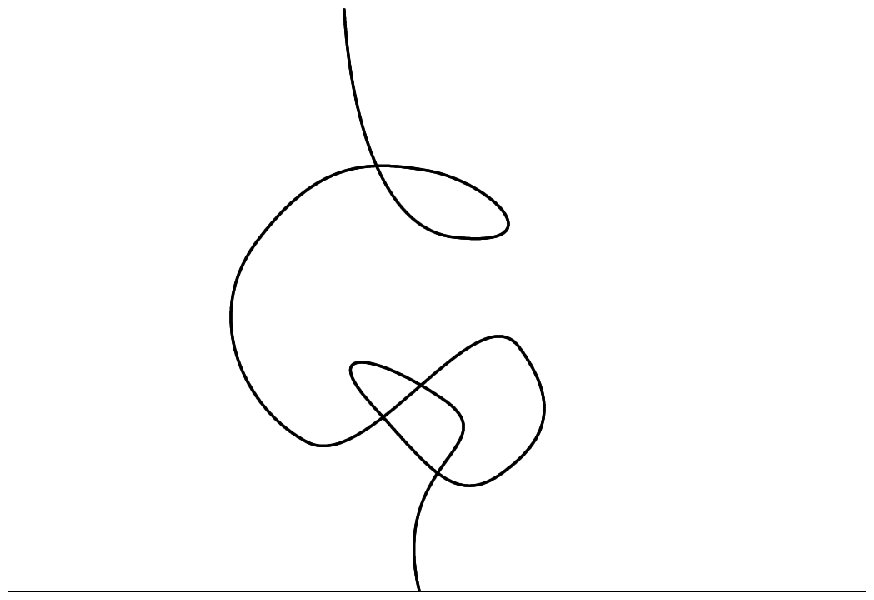}}
\caption{\label{f1.ex}
A curve} 
\end{figure}

\begin{figure}
\centerline{\includegraphics*[height=2in]{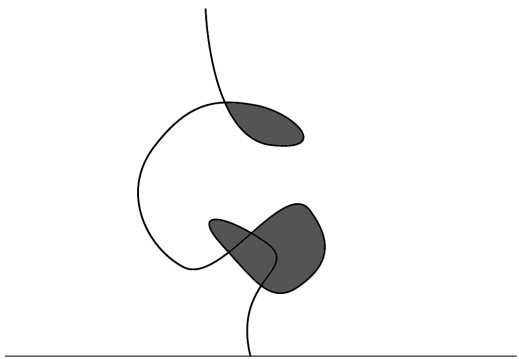}}
\caption{\label{f2.ex}
Its filling} 
\end{figure}

More generally, one can define fillings of other sets than curves: Fillings of the union of two curves, fillings of the union of fillings of curves etc.

It is easy to see that the law of the filling ${\cal F} (\gamma)$ of a 
random curve is characterized by the knowledge of the probabilities of the 
events ${\cal F} (\gamma) \subset H $. But (recall that the $H$'s are simply connected),
this is the same event as $\gamma \subset H$. In particular, we see that the 
law of the filling ${\cal F} (\gamma)$ of curves that satisfy conformal restriction 
is fully determined by the restriction exponent $\alpha$. 

This leads to the following definitions:
\begin {itemize}
\item
We say that a closed connected set $F$ connecting the origin to infinity in the upper half-plane
 is a filled if ${\cal F} (F) =F$. This is for instance the case if 
$F$ is the filling of a curve $\gamma$.
\item
We say that the random filled set $F$ satisfies conformal restriction, if its law is scale-inavriant (i.e. $\lambda F$
and $F$ have the same law) and if for all $H$, the 
law of $F$ given $F \subset H$ is identical to the law of $\Phi_H^{-1} (\H)$.
\end {itemize}
The same arguments as before show that if a random filled set $F$ satisfies conformal restriction, then there exists a positive exponent $\alpha$ such that for all $H$, 
\begin {equation}
\label {eq3}
P [ F \subset H] = \Phi_H' (0)^\alpha.
\end {equation}
And, conversely, for each positive $\alpha$, there exists at most one law of a random filled
set $F$ satisfying (\ref {eq3}) for all $H$. If it holds, then (and we leave this as an exercise), it follows that 
$F$ satisfies conformal restriction.
We call the law of $F$ the restriction measure with exponent $\alpha$ if it exists.
We know that when $\alpha=1$, this law exists: It is that of the 
filling of the conditioned Brownian motion.
 
Intuitively, the larger the exponent $\alpha$ is, the bigger the (possible) 
 corresponding random filled set $F$ should be (because then, $P[ F \subset H]$
 decreases with $\alpha$). For instance, suppose that $F_1$ and $F_2$ are two 
 independent filled sets with respective exponents $\alpha_1$ and $\alpha_2$.
 Then, define $F= {\cal F} (F_1 \cup F_2)$. Clearly,
 $F \subset H$ if and only if both $F_1$ and $F_2$ stay in $H$. Hence,
 $$
 P [ F \subset H] = \Phi_H' (0)^{\alpha_1} \Phi_H'(0)^{\alpha_2}
 $$ 
 and the random filled set $F$ satisfies conformal restriction with exponent 
 $\alpha_1 + \alpha_2$.
 
 This shows that for all positive integer $n$, the law of a random filled set satisfying conformal restriction with exponent $n$ exists: It can be constructed as  the filling of the union of $n$ independent Brownian motions.
 
 Also, suppose that there exists an exponent $\alpha_0$ such that the corresponding 
 law of a filled set satisfying conformal restriction exists, and is supported on 
 simple curves: Then, by filling unions of independent samples of this law and of 
 Brownian motions, one can construct the laws with exponents $n + m \alpha_0$ 
 for all integers $n\ge 0$ and $m \ge 1$. All these measures (except when $n=0$ and $m=1$
 of course) can not be supported on simple curves (if $n \ge 1$, it already contains 
 a Brownian motion with plenty of double points, if $m\ge 2$ it contains the union of two
 different simple curves). All this seems to indicate that when $\alpha$ is large, the 
 corresponding restriction measure is not supported on simple curves.
 
 On the other hand, it is not difficult to see that when $\alpha$ is too small, the 
 corresponding restriction measure does not exist. The reason is that the random set $F$ has to 
 connect the origin to infinity. In particular, it has to intersect the unit circle
 (or more precisely, the intersection of the unit circle and the upper half-plane).
 Hence, if $K_1= \{ \exp (i \theta) \ : \ \theta \in (0, \pi /2]\}$
 and $K_2 = \{\exp (i \theta) \ : \ \theta \in [\pi /2, \pi) \}$, then 
 $$
 P [ F \cap K_1 \not= \emptyset ] + P [F \cap K_2 \not= \emptyset ] \ge 1.
 $$
 Hence, 
 $$ 
 2 \Phi_{\H \setminus K_1}' (0)^\alpha \ge 1,
 $$
 which gives a lower bound to the admissible values of $\alpha$.
   
 We shall see that in fact, things are indeed as one might guess at this point:
 \begin {theorem}
 There exists a critical value $\alpha_0$ such that:
 \begin {itemize}
 \item
 If $\alpha < \alpha_0$, there is no random filled set satisfying conformal restriction with 
 exponent $\alpha$.
 \item
 There exists a random simple curve satisfying conformal restriction with exponent $\alpha_0$
 \item
 For all $\alpha > \alpha_0$, there exists a random filled set satisfying restriction 
 with exponent $\alpha$, and it is almost surely not a simple curve.
 \end {itemize}
 \end {theorem}
 This theorem partially generalizes Theorem \ref {main1} and its proof will proceed in several steps.
 As we shall see, the critical value is $\alpha_0= 5/8$ and corresponds to SLE$_{8/3}$. The scaling exponent $5/8$ 
 for the boundary behaviour of self-avoiding walks appeared first in the theoretical physics literature in a (slightly) different context in Cardy's paper \cite {Ca0}. 
 
It is worth stressing that this theorem gives a complete description of what the scaling limit of 
self-avoiding walks should be, assuming that it exists and is conformally invariant, because the formula 
$P [\gamma \subset H] = \Phi_H' (0)^{5/8}$ (for all $H$) gives the law of $\gamma$.
It is possible (but not trivial) to simulate very long self-avoiding walks in a 
half-plane, using a modified version of the 
pivot algorithm (basically, one has to find a Markov chains with the proper invariant measure, and let it 
run a sufficiently long time so that it reaches its stationary state). This procedure    
is explained in detail in \cite {Kenn1, Kenn2}.
This makes it possible to test numerically the conjecture that for very long rescaled 
self-avoiding walks in the half-plane, $P[\gamma \subset H]$ is close to 
$\Phi_H'(0)^{5/8}$ by a Monte-Carlo procedure. The results \cite {Kenn2} are very accurate, and to our knowledge, 
they are the most convincing evidence so far of the fact that the scaling limit of self-avoiding walks exist 
and are conformally invariant.
  
\section {The continuous intrinsic self-avoiding curve}

\subsection  {Introduction to SLE via SAW}

Suppose that the discrete measure on self-avoiding curves that was described in the introduction indeed  has a conformally invariant limit. We have already heuristically
argued that it should then satisfy conformal restriction.
 But, the discrete measure has an additional property that we shall now exploit. It is worthwhile stressing that this 
 additional property is not shared by the conditioned random walk/Brownian motion. 

If one knows the first $m$ steps of the walk, what is the law of the remaining steps?
In the case of the upper half-plane, one can heuristically argue as follows: The 
law of $\gamma$ is uniform. If one conditions the uniform measure, one obtains the 
uniform measure on the smaller set. In particular, the law of 
$S(m), S(m+1), \ldots$ will be the uniform measure on the self-avoiding walks from
$S(m)$ to infinity that stay in $\H \setminus \{S(0), \ldots,  S(m) \}$.
In the scaling limit, assuming its existence and conformal invariance, this property would become:

\medbreak
\noindent
{\sl $({\cal P})$:
Given $\gamma  [0,t]$, the conditional law of $\gamma [t, \infty)$ is identical to the 
law of $f^{-1} (\tilde \gamma)$, where $\tilde \gamma$ is an independent 
copy of $\gamma$ and $f$ is a conformal map from $\H \setminus \gamma [0,t]$ onto
$\H$ such that $f(\gamma_t) = 0$ and $f(\infty) =\infty$.}
 
\medbreak 
In the case of a finite domain $D$, the discrete property is even clearer. Recall that 
the mass of a walk $\omega$ of length $n$ from $A$ to $B$ in $D$ is proportional to 
$x^n$ for some well-chosen and lattice-dependent $x$. If one conditions on the 
first $m$ steps of $\omega$, then the mass of the future will be proportional to $x^{n-m}$
and the conditional measure is supported on the self-avoiding walks from $S(m)$ to $B$ in $D \setminus S[0,m]$.

Hence, one would like to find a random continuous curve $\gamma$ satisfying $({\cal P})$.
Also, the curve $\gamma$ should be symmetric with respect to the imaginary axis (the law of the image of $\gamma$ under this symmetry is identical to the law of $\gamma$). These two 
conditions are exactly those that did lead Oded Schramm \cite {S1} to define SLE, in the context
of loop-erased random walks (see \cite {S1}, or \cite {LLN,Lbook,Wstf} for a survey and introduction on 
SLE). Here is a very brief outline of how one constructs SLE building on this idea:

The first observation is that if $\gamma$ is a random simple curve from 0 to infinity in the 
upper half-plane, then it will be natural and useful to parametrize it in such a way that 
for each $t$, there exists a conformal map
$g_t$ from $\H \setminus \gamma [0,t]$ that  satisfies:
$$ g_t (z) = z + \frac {2t}{z}+ o (\frac 1z) $$
when $z \to \infty$. This can be thought of as a way to (re)-parametrize 
the curve in such a way that ``seen from infinity'' it grows at 
constant speed. Then, one can define $W_t = g_t (\gamma_t)$.
The previous property shows that (for the random curve that we are looking
for), the law of $g_t (\gamma [t, \infty))$ given $\gamma [0,t]$ is identical to the 
scaling limit of a self-avoiding walk from $W_t$ to infinity in $\H$. In other words, the
conditional law of $g_t (\gamma [t, \infty)) - W_t$ is identical to the law of $\gamma$.
 
When $s$ is very small, it is not difficult to see that for each fixed $z \in \H$
$$
g_s (z) = z + \frac {2s}z + o (s)
.$$ 
This is due to the fact that when $s$ is small, seen from $z$ and infinity, 
$\gamma [0,s]$ looks like a straight slit $[0, 2 i \sqrt {s}]$ at first order.
Recall that for a straight slit, one would have 
$$ g_s (z)= \sqrt {z^2 + 4 s }.
$$
Hence, this implies that (for general $\gamma$) at $t=0$, the time-derivative of $g_t(z)$ is $2/z$.

Similarly, when $t>0$ is fixed and $s$ is small, $g_t (\gamma [t,t+s])$
looks like a straight slit growing near $W_t$, and
$$g_{t+s} (z) = g_t(z) + \frac {2s} {g_t (z)-W_t} + o (s).$$
Hence,
\begin {equation}
\label {loewner}
\partial_t g_t (z) = \frac {2}{g_t (z) - W_t}.
\end {equation}

This equation is interesting, because it shows that it is in fact possible
to recover the curve $\gamma$ from the real-valued continuous 
function $W$: For each $z$, it suffices to solve the ordinary differential equation (\ref {loewner}) with $g_0 (z)= z$. This constructs the mappings $g_t$, and $\gamma$ then follows since $g_t^{-1} (\H) = 
\H \setminus \gamma [0,t]$. Hence, in order to construct a random simple curve $\gamma$,
it suffices to construct the corresponding random function $t \mapsto W_t$.

Property $({\cal P})$ implies exactly that given $W[0,t]$, the law of 
$(W(t+s)-W(t), s \ge 0)$ is identical to the law of an independent copy of $W$.
In other words, $W$ is a continuous Markov process with independent increments.
Symmetry shows that $W$ and $-W$ have the same law. Hence, the only possibility is
that $W$ is real Brownian motion. More precisely, there exists a constant variance
$\kappa \ge 0$ such that $W_t / \sqrt {\kappa}$
is a standard Brownian motion $B$. 

To sum up things, we have just seen that if the (simple)
scaling limit of the self-avoiding curves
exist and are conformally invariant, then they can be constructed as follows, for some 
given constant $\kappa$:
Define $W_t = B_{\kappa t}$ where $B$ is ordinary real-valued Brownian motion. Then, solve
for each fixed $z$ the equation (\ref {loewner}) with initial data $g_0 (z) = z$.
This defines the mappings $t \mapsto g_t (z)$. Since one can do this for each $z$, this
procedure defines the conformal maps $z \mapsto g_t (z)$ for each fixed $t$. Then,
$\gamma$ is constructed by 
$\gamma(0,t]= \H \setminus g_t^{-1} (\H)$, or more precisely by
$$ 
\gamma (t) = g_t^{-1} ( W_t)
$$
if $g_t^{-1}$ extends continuously to $W_t$ (and it turns out to be almost surely the case).
The curve $\gamma$ is called the Schramm-Loewner evolution (SLE) with parameter $\kappa$.
Actually, it is called chordal SLE to indicate that it goes from one point of the boundary to another 
boundary point of the domain ($\H$ here), as opposed to other versions (radial, whole-plane), but since chordal 
SLE will be the only one that we will study in these notes, we just call it SLE.

For a general definition/introduction to SLE, see for instance \cite {Wstf, Lbook}, or the original paper \cite {S1}.
The following results can be proved, but they are not easy. The purpose of these lectures is not 
to focus on them, so we just list them, without further justification (see \cite {RS} for the original proofs):

\begin {itemize}
\item
For all $\kappa \le 4$, this procedure does indeed (almost surely) construct a simple curve $\gamma$
(\cite {RS}).
\item
For all $\kappa > 4$, this procedure does (almost surely) construct a continuous curve $\gamma$, but this 
curve is not simple. It has double points.
\item
For all $\kappa \ge 8$, this procedure does almost surely construct a continuous space-filling curve $\gamma$.
\end {itemize}
Also, it can be shown (except in the case where  $\kappa=4$, which is still open), that the Hausdorff dimension
of $\gamma$ is almost surely $1 + \kappa /8$ when $\kappa \le 8$.  See \cite {Bef} for the general case.

\subsection {Conformal restriction for SLE$_{8/3}$}

We are now going to combine this with conformal restriction, that the random simple curve 
$\gamma$ (that would be the scaling limit of the infinite half-plane self-avoiding walk) should also satisfy:
Let us fix an $H$. We know that for some exponent $\alpha$, 
$$
P [ \gamma \subset H ] = \Phi_H' (0)^\alpha.
$$
If we know $\gamma [0,t]$, what is the conditional probability that $\gamma \subset H$ when $\gamma$ is 
an SLE?
Of course, this question is non-trivial only if $\gamma [0,t] \subset H$.
Let us map the future of $\gamma$ by the uniformizing map $g_t$. The conditional law
of $g_t (\gamma [t, \infty)) - W_t$ is the same as the law of $\gamma$.
In particular, for any $H'$,
$$
 P [  g_t (\gamma [t, \infty)) - W_t  \subset H' \mid \gamma [0,t] ] 
= \Phi_{H'}'(0)^\alpha.
$$
Note that $\gamma [t,\infty) \subset H$ means that 
$ g_t (\gamma [t, \infty)) \subset g_t (H)$. 
Hence, 
$$
P [ \gamma \subset H \mid \gamma [0,t] ] = \Phi_{g_t (H) -W_t}' (0)^\alpha
=\Phi_{g_t(H)}'(W_t)^\alpha.$$
This means that this last quantity must be a martingale, and 
this has to hold for any $H$.

It so happens that this only holds for one specific choice of $\kappa$. Namely
$\kappa=8/3$ and the corresponding value of $\alpha$ is then $5/8$.
Let us now outline how the computation goes.
It can be performed directly and rigorously 
for macroscopic $\H \setminus H$ (we will mention this proof later on), but at least on the 
heuristic level, one can also focus on the case where 
$H = \H \setminus [x, x+i \delta]$ for an infinitesimally small $\delta$.
Recall that the conformal map $\Phi_H$ in this case is
$$
\Phi_H  (z) = \sqrt { (z-x)^2 + \delta^2 } - \sqrt {x^2 + \delta^2}.
$$
The probability to avoid the infinitesimal slit is therefore roughly 
$$ 
\Phi_H'(0)^\alpha = (1 + \frac {\delta^2}{x^2})^{- \alpha/2}
= 1 - \frac {\alpha}{2x^2} \delta^2 + o (\delta^2)
$$
when $\delta \to 0$.

After a small time $t$, let us look what the conditional probability to 
hit this infinitesimal slit becomes: 
After mapping by the map $g_t$, it is the probability that an SLE started from
$W_t$ hits $g_t ( [x, x + i \delta])$ that  (at first order) is the same as the 
infinitesimal slit 
$[g_t (x), g_t(x)+ i \delta g_t' (x)]$. The value of this quantity is 
(by the same computation as before)
$$
 1 - \frac {\alpha g_t'(x)^2}{2(g_t(x) - W_t)^2} + o (\delta^2).
$$ 
In other words, we need 
$$ g_t'(x)^2 / (g_t(x) - W_t)^2$$
to be a (local) martingale. Recall that $W_t = \sqrt {\kappa} B_{ t}$,
that $g_t (x)$ and $g_t' (x)$ (by formal differentiation with respect to $x$) 
are $C^1$ functions of time:
$$
\partial_t g_t (x)= \frac 2 {g_t(x) - W_t}
\hbox { and }
\partial_t g_t' (x) = \frac {-2 g_t'(x)} {(g_t (x) - W_t)^2}.
$$
It\^o's formula (loosely speaking, Taylor's expansion using the 
fact that the mean value of $W_t^2$ is $\kappa t$) 
shows that the drift term of $g_t'(x)^2 / (g_t (x) - W_t ) ^2$
(i.e. the mean first term in $t$ when $t \to 0$) is 
\begin {eqnarray*}
&&\frac {-4 g_t' (x)^2}{ (g_t(x) - W_t)^4}
 +
 \frac {-4 g_t' (x)^2}{ (g_t(x) - W_t)^4}
 + \frac {\kappa}{2} \times \frac { 2 \times 3 \times g_t'(x)^2}{(g_t (x) - W_t)^4 }\\
 && = 
 (-8 + 3 \kappa) \frac {g_t' (x)^2}{ (g_t(x) - W_t)^4}
\end {eqnarray*}
that vanishes only when $\kappa = 8/3$. 

The value of $\alpha$ can then determined by inspection of the higher-order terms. This infinitesimal 
approach was shown in \cite {FW} to be related to highest-weight representations of the algebra of polynomial 
vector fields on the unit circle (these representations were those used in the theoretical physics literature to predict the value of critical exponents). We will come back to this later.

\section {One-sided restriction}

\subsection {Definition} 
Suppose now that a random curve (or filled set) $F$ from the origin to infinity in the 
upper half-plane $\H$ satisfies the following weaker form of conformal restriction:
For any simply connected $H_+ \subset \H$ such that $\H \setminus H_+$ is bounded, and 
bounded away from the whole negative half-line, the conditional law of $\Phi_{H_+} (F)$
given $F \subset H_+$ is identical to the law of $F$.
We then say that $F$ satisfies one-sided restriction.

The difference is that we impose that $H_+$ has the negative half-line on its boundary.
In doing that, we break the $\sigma$-symmetry with respect to the imaginary line. 
In the sequel, when we use the notation $H_+$, we will always implicitly mean for such sets.

Of course, if $F$ satisfies restriction then it satisfies also one-sided restriction.
Conversely, if $F$ satisfies one-sided restriction, then, the same arguments as in the 
two-sided case end up showing that there exists $\alpha>0$ such that for all $H_+$,
\begin {equation}
\label {ex2}
P [ F \subset H_+ ] = \Phi_{H_+}' (0)^\alpha.
\end {equation}
As before, this relation does not fully characterize the law of $F$. It does characterize the 
law of ``its right-boundary''. This can be defined in terms of its one-sided filling, i.e.
the set of points in $\H$ that are separated from $\R_+$ by $F$. The boundary of this 
one-sided filling consists of $\R_-$ and of a curve $\gamma_+ (F)$, that we call the right-boundary of $F$.
It is easy to see that the relation (\ref {ex2}) characterizes the law of $\gamma_+ (F)$.
Note that it is a priori not clear that $\gamma_+$ is a curve, even less that it is a simple curve, but it 
will turn out that it is indeed the case. We call it (if it exists) a one-sided restriction curve with
exponent $\alpha$.

\begin{figure}
\centerline{\includegraphics*[height=2in]{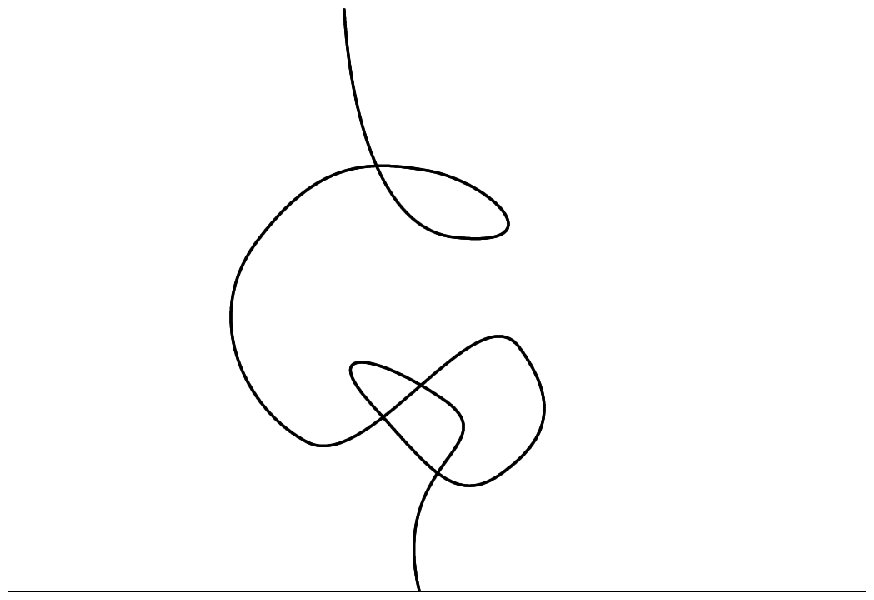}}
\caption{\label{f3.ex}A curve} 
\end{figure}

\begin{figure}
\centerline{\includegraphics*[height=2in]{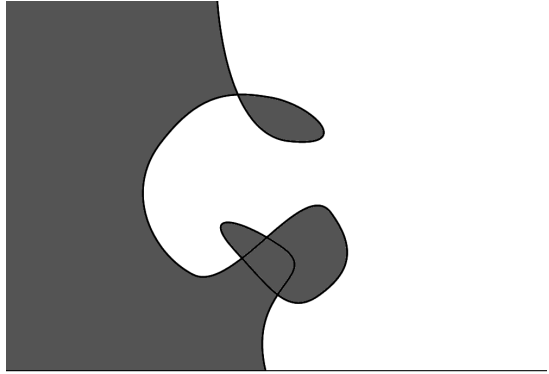}}
\caption{\label{f4.ex}
Its one-sided filling} 
\end{figure}

We shall see that for all $\alpha \ge 5/8$, the two-sided restriction measure with exponent 
$\alpha$ exists. Hence, the one-sided restriction curve with the same exponent exists  too (just take the right-boundary of the sample of the two-sided measure).
Recall that the non-existence of the (two-sided) restriction measure for small $\alpha$ was based on a symmetry 
argument. This does not apply to the one-sided case. In fact, we shall see that the one-sided 
restriction curve exist for all $\alpha$ and we 
will show four very different ways to construct these one-sided restriction measures.

\subsection {Reflected Brownian excursions}

We start with the example of the Brownian motion, that it started from the origin, 
reflected orthogonally on the negative half-axis (so that it stays in $\H$) and conditioned 
never to hit the positive half-line (this will, as before, imply that the motion is transient and
tends to infinity).
It is in fact convenient to use a reflection and conformal invariance argument to construct 
this
process starting from a sample $B$ of $P_{\H,0,\infty}^{BM}$.
Modulo time-change, the law of the path $B^2$ (the path of the square of the complex-valued process $B$) is 
$P_{\C \setminus \R_+, 0, \infty}^{BM}$. We now define $Z$ by reflection of 
$B^2$ with respect to the real axis: The real part of $Z$ is the real part of $B^2$, and the imaginary part 
of $Z$ is the absolute value of the imaginary part of $B^2$.
The path of $Z$ is (modulo time-change) that of a Brownian motion in the upper half-plane, that is 
orthogonally reflected on the negative half-line, and conditioned not to hit the positive
half-line. 

Because $B$ satisfies restriction, it is easy to see that $Z$ satisfies one-sided restriction.
For each $H_{+}$, conditioning $Z$ to stay in $H_{+}$ means to condition $B$ to stay in some
set $H$ (and does not change the law, modulo conformal invariance). Furthermore, it is easy (and left to the 
reader) to see that the obtained one-sided exponent is $1/2$. More generally, this procedure can be applied to 
any two-sided restriction measure, and produces the one-sided restriction measure with half its exponent.

This has the following rather surprising consequence: The right-boundary of the union of two copies 
of $Z$ (conditioned reflected Brownian motion) has the same law as the right-boundary of $B$ (conditioned
Brownian motion). We shall see plenty of such identities in law between the right-boundaries of sets that 
are constructed in very different ways.

It is interesting to note that the process $Z$ is a scale-invariant Markov process. Reflected Brownian motion 
is conformally invariant, so that, just as for the ordinary Brownian motion, one can a priori argue that this
conditioned reflected Brownian motion does satisfy one-sided conformal restriction. The same argument applies 
also even when the reflection is not orthogonal. We will not describe non-orthogonal reflection in detail
here (see e.g. \cite {VW}) but we mention that for each $\theta \in (0, \pi)$, there exists a Brownian motion
in the upper half-plane that is reflected on the real line with an angle $\theta$ (i.e. the push when it hits the 
real line is 
proportional to $\exp(i \theta)$). If one starts this process from the origin and (appropriately) conditions it never to hit the positive real axis, then one obtains a random path from the origin to 
infinity that hits the negative half-line, but not the positive half-line. It is then possible (and not difficult)
to see that:

\begin {theorem}
The right-boundary of reflected Brownian motion (with angle $\theta$ on the negative half-axis) conditioned to 
never hit the positive half-axis satisfies one-sided restriction with exponent $\alpha =1  - \theta/ \pi$.
\end {theorem}

When $\theta$ is close to $0$, then the reflection pushes the motion very strongly towards the origin (and therefore 
towards the positive half-line, which the motion tries to avoid because of the conditioning), so that in the 
limit where $\theta=0$, one obtains simply Brownian motion conditioned to avoid the whole real line and the 
restriction measure with exponent $1$.
Similarly, it is not difficult to see that for all $H_+$, in the limit where $\theta \to \pi$, 
the probability that $B$ stays on $H_+$ goes to one (and the limiting exponent is $0$).

Note that the right-boundary of the union of independent sets satisfying one-sided restriction does also satisfy one-sided restriction, and (as in the two-sided case), the corresponding exponent is the sum of the exponents.
In particular, for all $\alpha >0$, it is possible to find $\alpha_1, \ldots, \alpha_n$ in $(0,1)$ 
such that $\alpha_1 + \cdots + \alpha_n= \alpha$, and therefore to construct a set satisfying one-sided 
conformal restriction with exponent $\alpha$, as the right-boundary of the union of $n$ conditioned 
reflected Brownian motions. This also implies that this right-boundary is almost surely a path (because 
planar Brownian motion can be proved to have no double cut-points, see  \cite {BL}).

\begin {corollary} 
\label {C1}
For all $\alpha >0$, the one-sided restriction path with exponent $\alpha$ does exist and it is a simple path.
\end {corollary}
    
We can now already use this to prove the following result on two-sided restriction measures:

\begin {corollary}
\label {C2}
For all $\alpha < 5/8$, the two-sided restriction measure with exponent $\alpha$ does not exist.
\end {corollary}

\noindent
{\bf Proof.}
Suppose that $\gamma_+$ is the right boundary of a two-sided restriction 
measure with exponent $\alpha < 5/8$. By symmetry, the probability that it passes to the 
right of $i$ (i.e. that it separates $i$ from the positive half-line in $\H$) is at least $1/2$. 
By adding an independent conditioned reflected Brownian motion (with appropriately chosen 
angle) and taking the right-boundary of the union, one obtains the right-boundary of a
one-sided restriction sample with exponent $5/8$, and the probability that it
passes ``to the right of $i$'' is then strictly larger than $1/2$. But we know that the right-boundary 
of the two-sided restriction measure with exponent $5/8$ can be realized as an SLE$_{8/3}$ curve itself 
(we know it is a simple curve, so it is its own right boundary), and by symmetry, it has a probability 
$1/2$ to pass to the right of $i$. 
This leads to a contradiction.

\subsection {Poisson clouds of Brownian excursions}

We now describe another way to construct the one-sided restriction measures with exponent $\alpha$.
It is in fact related to the previous one, and corresponds to the limiting case where one decomposes 
$\alpha$ in $n$ times $\alpha/n$, so that the restriction sample is constructed as the union of a lot of independent 
conditioned reflected Brownian motions with a very steep angle (here, close to $\pi$). 

For each real $x$, there is a natural infinite measure on Brownian paths started from $x$ in the upper half-plane:
$$
\mu_{\H,x} = \int_\R \frac {dy}{y^2} P_{\H, x,x+y}^{BM}. 
$$
This is simply the rescaled limit when $\eps \to 0$, of the law of Brownian motion started from 
$x+i\eps$, and killed at its first exit of $\H$ (the previous integral corresponds simply to a decomposition 
according to the location of the exit point $y$).
 
We now construct a measure $\mu_\H$ with non-prescribed starting point by integrating the starting 
point according to the Lebesgue measure on $\R$: $\mu_\H = \int_\R dx \mu_{\H,x}$.
This is the Brownian excursion measure, as constructed in \cite {LW2}.
Finally, we restrict this measure to the set of paths that start and end on the negative half-line. We call $\mu^-$
the obtained measure. In other words:
$$
\mu^- = \int_{\R-} dx \int_{-\infty}^{-x} \frac {dy}{y^2} P_{\H, x,x+y}.
$$

Suppose now that the set $H_+$ is given, and define $\mu^-_{H_+}$ to be $\mu^-$
restricted to those paths that stay in $H_+$.
Then, one can take the image of this measure under $\Phi_{H_+}$. It is a simple exercise 
(using the fact that $P_{\H,x,x+y}$ satisfy restriction with exponent $1$), to see that this 
image measure is identical to $\mu^-$ itself. The point is that the $\Phi'(x) dx$ term due to the 
map $\Phi$ is balanced out exactly by conformal restriction.
 The measure $\mu^-$ therefore satisfies some generalization of one-sided conformal
restriction (it is not exactly the same property as before though, because it is an infinite measure, so that conditioning does not make sense anymore). 

For each (even infinite) measure $\mu$ on a state-space $S$, one can construct a Poissonian realization of $\mu$.
This is a random countable (or finite) collection $(x_j, j \in J)$ of elements of $S$. Its law is characterized by 
the following two facts:
\begin {itemize}
\item
For any disjoint measurable subsets $A$ and $A'$ of $S$, the events 
$\{ \exists j \in J \ : \ x_j \in A \}$ 
and 
$\{\exists j \in J \ : \ x_j \in A' \}$ are independent.
\item
For any measurable $A$, the number of $j$'s such that $x_j \in A$ 
is a Poisson random variable with mean $\mu (A)$.
\end {itemize}
Hence, for each $\beta >0$, one can define a random Poissonian 
realization $(\gamma_j)$ of $\beta \mu^-$. Since the measure $\mu^-$ is infinite, this 
is an random infinite collection of Brownian curves that start and end on the 
negative half-line. For each disjoint compact intervals $I$ and $I'$, the number of 
curves that start in $I$ and end in $I'$ is almost surely finite. But the 
number of curves that start and end in $I$ is almost surely infinite. There 
are only finitely many ``macroscopic'' curves, and infinitely many microscopic ones.
We call $Q^\beta$ the law of this random collection.

Not that the definition of the Poissonian realization implies immediately that the 
union of two independent realizations of $Q^\beta$ and $Q^{\beta'}$
is a realization of $Q^{ \beta + \beta'}$.

Suppose now that $(\beta_j, j \in J)$ is a realization of $Q^\beta$. One can consider the 
right-boundary $\gamma_+$ of $\cup_j \gamma_j$. 

\begin {theorem}
The right-boundary $\gamma_+$ of $\cup_j \gamma_j$ satisfies one-sided conformal restriction 
with exponent $c\beta$ (for some given constant $c$).
\end {theorem}

\noindent
{\bf Proof.}
It suffices in fact to check that $\gamma_+$ satisfies conformal restriction (the relation between 
the exponent and $\beta$ then follows from the above-mentioned additivity property), and this follows 
rather readily from the ``conformal restriction'' property of $\mu^-$.

\medbreak
Note that Corolaries \ref{C1} and \ref{C2} can also be deduced from this alternative construction of 
one-sided restriction measures.

\subsection {Some remarks}

We have now seen three different but equivalent constructions of this SLE$_{8/3}$ curve:
\begin {itemize}
\item The SLE construction, via Loewner's equation.
\item The right boundary of a reflected Brownian motion with angle $\theta=3 \pi /8$ on the negative 
half-line, conditioned not to intersect the positive half-line.
\item The right-boundary of a Poisson cloud of Brownian excursions attached to the negative half-line.
\end {itemize}

While in the first case, it is clear that this produces a random object that is symmetric with respect to 
the imaginary half-axis, it can seem quite amazing that the two latter constructions do. In particular, this
shows for instance that in the second construction if one sees only the outer boundary, one cannot tell whether 
it has been generated as the right boundary of a Brownian motion reflected on the negative half-axis, or as the 
left boundary of a Brownian motion reflected on the positive half-axis.
This shows that ``the outer boundary of planar Brownian motion is locally symmetric'': If one only sees a piece 
of this boundary, one cannot tell on which side the Brownian motion is.
A similar observation follows from the fact (that only uses two-sided restriction) that the filling of the 
union of eight independent SLE$_{8/3}$ has the same law as the filling of the union of 5 independent conditioned
Brownian motions (they both satisfy conformal restriction with exponent $5$). 

In any case, this shows that:
\begin {corollary}
If the random simple curve $\gamma$ satisfies one-sided restriction, then its Hausdorff dimension is almost surely $4/3$.
If a random filled set satisfies two-sided restriction, then the Hausdorff dimension of its outer boundary 
is almost surely $4/3$.
\end {corollary}

Let us also note that the construction of one-sided restriction measures via Poisson clouds of excursions show that
it is possible to construct on the same probability space an increasing family $(F_\alpha)_{\alpha >0}$ of right-filled sets such that for each fixed $\alpha$, $F_\alpha$ satisfies one-sided restriction with exponent $\alpha$
(the point is that it is an increasing family). 
In the two-sided case, (to my knowledge) the corresponding problem is open:

\begin {question}
Is it possible 
to define an INCREASING family of filled sets $(F_\alpha)_{\alpha \ge 5/8}$ satisfying two-sided restriction with 
respective exponents $\alpha$?
\end {question}

\section {Restriction defect when $\kappa \not= 8/3$}

\subsection {The martingale}
As we have already seen, when $\kappa \not= 8/3$, SLE does not satisfy conformal restriction.
Let us start with analytic considerations, and let us see what goes wrong with the ``global'' 
proof of conformal restriction as soon as $\kappa \not= 8/3$.

Suppose that $H$ is fixed, and that $\gamma$ is an SLE$_\kappa$. Define the maps $g_t$ as before.
Also, define
$$ 
h_t = \Phi_{g_t (H)}.
$$
This means that $h_t \circ g_t$ is a conformal map from $H \setminus \gamma [0,t]$ onto $\H$. 
When $\kappa = 8/3$, let us look at the conditional probability that $\gamma \subset 
H$ given $\gamma [0,t]$. When $\gamma [0,t] \subset H$, this means that 
$g_t (\gamma [t, \infty)) \subset g_t (H)$. In other words, this conditional probability is 
equal to $h_t' (W_t)^\alpha$, where $\alpha= 5/8$.

The actual proof of this fact proceeds as follows. First, one has to note that
the mapping $(t,z) \mapsto h_t' (z)$ is $C^1$ with respect to $t$ and analytic with respect to $z$.
One works out the derivative with respect to time of 
$h_t'(z)$ in terms of the mapping $h_t$. A simple computation based (almost) only on the expression of the 
derivative of the composition of maps shows that
$$
\partial_t h_t (z) =  \frac { 2 h_t' (W_t)^2 } {h_t (z) - h_t (W_t)}
- \frac {2 h_t' (z)} {z - W_t}.
$$
One can formally differentiate with repect to $z$, and get the expression for
$\partial_t h_t'(z)$.
In the limit where $z \to W_t$, we then get that 
$$
(\partial_t h_t) (W_t) = -3 h_t''(W_t) \hbox { and } (\partial_t h_t') (W_t) = 
\frac {h_t'' (W_t)^2} {2 h_t' (W_t)} - \frac { 4 h_t''' (W_t)} {3}.
$$
If one looks at 
$ d ( h_t' (W_t))$, one has to use a 
 slightly extended version of It\^o's formula (i.e. Taylor's expansion using the fact that $(dW_t)^2 = \kappa dt$):
$$ d (h_t' (W_t)) = h_t'' (W_t) dW_t + ( \partial_t h_t' ) (W_t) dt + \frac {\kappa}2 h_t''' (W_t) dt.
$$ 
Similarly, if one looks at the variation of $h_t'(W_t)^\alpha$, one gets that
\begin {eqnarray*}
d (h_t' (W_t)^\alpha)
&=&
\alpha h_t' (W_t)^{\alpha -1} h_t''(W_t) dW_t 
\\
&&+ \Big(\frac {\alpha (\alpha-1)}2 h_t' (W_t)^{\alpha-2} h_t''(W_t)^2 + 
\alpha h_t'(W_t)^{\alpha - 1} h_t'''(W_t) \Big) (dW_t)^2\\
&&+ \alpha h_t' (W_t)^{\alpha-1} \partial_t h_t' (W_t) dt.
\end {eqnarray*}
Hence, one ends up with the fact that
\begin {eqnarray*}
\frac {d(h_t' (W_t))^\alpha}{ \alpha h_t'(W_t)^\alpha} 
&=&\frac { h_t''(W_t)}{h_t'(W_t)}  dW_t 
\\
&& + \Big( \frac { (\alpha-1) \kappa + 1} 2 \frac {h_t''(W_t)^2}{ h_t'(W_t)^{2}}
\\
&&
+ ( \frac \kappa 2 - \frac 4 3 ) \frac {h_t'''(W_t)}{ h_t'(W_t)} \Big) dt.
\end {eqnarray*}
In the special case where $\alpha= 5/8$ and $\kappa =8/3$, this is indeed a (local) martingale, because the 
$dt$ term vanishes. This gives the direct analytical proof of the fact that $SLE_{8/3}$ satisfies conformal 
restriction with exponent $5/8$.

When $\kappa \not= 8/3$, there is no choice of $\alpha$ and $\kappa$ that turns 
$h_t'(W_t)^\alpha$ into a local martingale. But it is natural to define
$$\alpha (\kappa) =  \frac {6 - \kappa }{ 2 \kappa}.$$
For this value,
\begin {equation}
\label {Sch}
d (h_t' (W_t)^\alpha)= \alpha h_t' (W_t)^\alpha 
\left( \frac {h_t''(W_t)}{h_t'(W_t)} dW_t + (\frac \kappa 2 - \frac 4 3) S_{h_t}(W_t) dt \right) .
\end {equation}
Here and in the sequel, $S_f$ stands for the Schwarzian derivative of a map $f$ defined by
$$
S_f = \frac {f'''}{f'} - \frac {3 (f'')^2} {2 (f')^2}.
$$
This is a well-known quantity in the theory of conformal maps. In our context, we will 
only use it for conformal maps $h_t$ and at a boundary point. 
Recall that $h_t$ is a conformal map (normalized) at infinity that removes $\H \setminus g_t (H)$
(i.e. it is a map from $g_t (H)$ onto $\H$). One can think of $-S_{h_t} (z)$ when $z$ is on the real line, as
a (conformally invariant) way to measure the size of this removed set, seen from $z$ (this will be clear with
the Brownian loops interpretation that we will give a little later). In particular, $S_{h_t} (W_t)$ is a negative quantity.
 
Equation (\ref {Sch}) shows immediately that 
\begin {equation}
\label {themartingale}
M_t = h_t'(W_t)^\alpha \exp ( \frac  \lambda 6 \int_0^t S_{h_s} (W_s) ds )
\end {equation}
is a local martingale where 
$$\lambda (\kappa) = (8 - 3 \kappa) \alpha = \frac {(8-3 \kappa)(6- \kappa)}{2 \kappa}.$$
This is simply due to the fact that 
$$
dM_t = d(h_t'(W_t)^\alpha) \exp ( \frac  \lambda 6 \int_0^t S_{h_s} (W_s) ds )
+ \frac \lambda 6 S_{h_t} (W_t) M_t dt.
$$

Since $S_{h_t} (W_t)$ is negative, the local martingale $M$ is positive and bounded by one if $\kappa \in [0, 8/3]$. 
This implies that it converges in $L^1$, and therefore that
 $M_t$ is the conditional expectation of some random variable $M_\infty$ given $\gamma[0,t]$.

\subsection {The $\kappa=2$ case}

It is useful to focus on the case where $\kappa=2$ (note that in this case, the martingale is bounded).
SLE$_2$ is one of the few special cases where it has now been rigorously established that it is the scaling limit of a discrete model from statistical physics: SLE$_2$ is the scaling limit of loop-erased random walk. More precisely, the 
proofs in \cite {LSWlesl} can be adapted to show that if one erases chronologically the loops of the random walk in a lattice (say, the square lattice) conditioned to remain forever in the upper half-plane (the path that converges to the
conditioned Brownian motion) and considers the scaling limit of this loop-erased random walk, one obtains chordal 
SLE$_2$.
The loop-erasing procedure is chronological in the sense that if $(Z_n, n \ge 0)$ is the conditioned random 
walk, then the loop-erased path $(L_p, p \ge 0)$ is defined by
$L_0= 0$, 
and for each $p\ge 0$, $n_p = \sup \{ n \ : \ Z_n = L_p \}$, and 
$$ L_{p+1} = Z_{n_p +1}$$
(all the loops from $L_p$ back to itself have been erased).

In the discrete setting, it is therefore possible to couple a loop-erased random walk together with the 
(conditioned) random that was used to construct it. It is in fact not difficult to understand the conditional 
law of $Z$ if one knows $L$: One has to add random loops back on the top of $L$ in some appropriate way.

\begin{figure}
\centerline{\includegraphics*[height=3in]{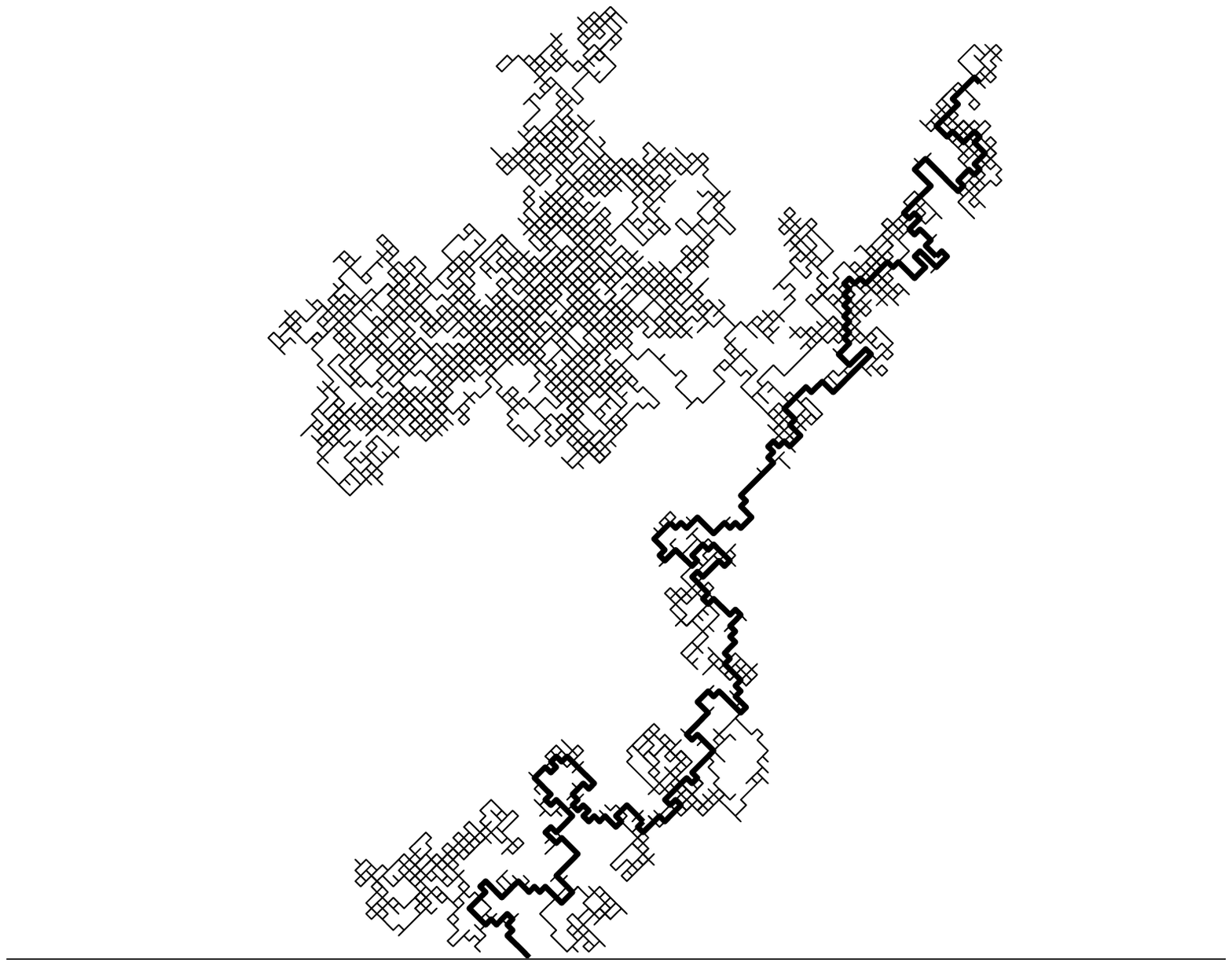}}
\caption{\label{f.le}
Beginning of a simple random walk and its loop-erasure}
\end{figure}

In the scaling limit, this should remain valid: It is possible to couple an SLE$_2$ with  
a conditioned Brownian motion, in such a way that in some sense, the SLE$_2$ is the loop erasure of 
the Brownian motion (here, we avoid some subtle open questions such as: Is this ``loop-erasure'' deterministic
in the scaling limit?). Conversely, if one adds the Brownian loops back on top of the SLE$_2$ curves, then one 
constructs exactly the conditioned Brownian motion. It is possible to describe fully and rigorously this procedure 
to put Brownian loops back on top of a curve. This construction may seem somewhat complicated at first sight, but it is 
in fact the direct generalization of the discrete case, that makes use of conformal invariance.

One starts with an infinite measure on paths $\nu$ from $0$ to $0$ in the 
upper half-space. It is a renormalized version of $P_{\H, 0, 0}^{BM}$.
Then, one considers the measure $dt \otimes \nu$ where $dt$ is the Lebesgue measure on $\R_+$,
and a Poissonian realization $((t_j, \eta_j), j \in J)$ of $dt \otimes \nu$.
One can view this realization as a random family of Brownian loops $(\eta^t, t\ge 0)$
such that for all but countably many $t$'s, $\eta^t=\emptyset$, whereas $\eta^t = \eta_j$
when $t = t_j$ for $j\in J$.
On each finite time-interval, there are countably many loops, but only finitely many macroscopic ones (of diameter larger than one, say).
 
Then, if $(\gamma_t, t\ge 0)$ is the SLE$_2$, and $g_t$ the corresponding conformal maps, one adds to the SLE the loops $l_{t_j}= g_{t_j}^{-1} (W_{t_j} + \eta_j)$. One transforms conformally the loop $\eta_j$ (that is a loop from $0$ to $0$ in $\H$) into a loop from $\gamma_{t_j}$ to $\gamma_{t_j}$ in $\H \setminus \gamma [0, t_j]$.
 
Since (SLE$_2$ + the loops) should form a conditioned Brownian motion, 
the obtained path satisfies conformal restriction with exponent $1$.
In particular, 
the probability that neither the SLE nor the loops exit a given $H$ is $\Phi_H'(0)$. 
Let $E_H$ denote this event. 
Let us now try to understand the conditional probability of $E_H$ given $\gamma [0,t]$.
Note that we do not have the knowledge of the loops (even for $t_j \le t$) here.
There are two contributions:

First, it means that none of the loops that have been added before time $t$ does exit $H$. The probability that one 
adds a loop in the time-interval $[t, t+dt]$ that exits $H$ is $- dt S_{h_t} (W_t) /3$
(this is not really surprising, $-S_{h_t} (W_t)$ is a quantity that measures the ``size'' of $h_t^{-1} (\H)$ seen from $W_t$ in $\H$). Hence, because of the Poissonian procedure, 
the probability that none of the loops that have been added before $t$ exits $H$ is simply
$$ \exp ( \int_0^t S_{h_s} (W_s) ds/3 ). $$  

Second, there is the probability that in the future, the SLE + loops do stay in $H$. This is exactly the probability 
that the Brownian motion started from $\gamma_t$ and conditioned to stay forever in $\H \setminus 
\gamma [0,t]$ stays also in $H$ (recall the definition of loop-erased random walk). 
By conformal invariance, this is exactly $h_t' (W_t)$.

These two contributions are conditionally independent (because of the Poissonian procedure, and the Markov property of the conditioned Brownian motion), so that in the end:
$$
P [ E_H \mid \gamma [0,t] ] 
= h_t' (W_t) \exp ( \int_0^t S_{h_s} (W_s) ds /3 ) .
$$
This is exactly the martingale $M_t$ when $\kappa =2$.
When $t \to \infty$, the limit is
$P [E_H \mid \gamma [0,\infty]]$. When $\gamma$ exits $H$, this quantity is zero. When 
$\gamma \subset H$, then it is the probability that no added loop to this specific $\gamma$ exits
$H$.

\subsection {The general case}
 
 When $\kappa \in (0, 8/3)$, one can interpret the martingale $M_t$ in a very similar way to the 
 $\kappa=2$ case. One adds to SLE$_\kappa$ a certain Poissonian cloud of Brownian loops exactly as 
 before, except that one changes the value of the parameter $\lambda$ (doubling $\lambda$ means for instance to take 
 the union of two independent clouds corresponding to $\lambda$) into $(6- \kappa)(8- 3 \kappa)/2 \kappa$.

Then, and this can seem surprising, the obtained path satisfies conformal restriction. The  restriction exponent
is
$\alpha = (6- \kappa)/2 \kappa$. The martingale $M_t$ can then still be interpreted as 
$$ 
M_t = P [ E_H \mid \gamma [0,t]]
,$$
where $E_H$ is the event that (SLE+loops) stays in $H$.

When $\kappa$ is close to $8/3$, the density $\lambda$ of the cloud vanishes, and the exponent is close to $5/8$
(from above), which is not surprising since when $\kappa = 8/3$, the SLE itself (without extra loops) satisfies 
conformal restriction. When $\kappa$ gets smaller, one has to add more loops, and the exponent of the obtained 
restriction measure gets larger. In the limit where $\kappa$ is very small, then $\lambda$ and $\alpha$ both go infinity. This is not surprising as when $\kappa$ is very small the SLE is very close to a straight line, and 
therefore very far from a set satisfying conformal restriction.

Note that for each $\alpha \ge 5/8$, there exists a value $\kappa \le 8/3$ such that $\alpha (\kappa) = \alpha$.
The conformal restriction measure with exponent $\alpha$ can therefore be constructed by this loop-adding procedure
and proves that:

\begin {theorem}
For each $\alpha \ge 5/8$, the (two-sided) conformal restriction measure with exponent $\alpha$ exists.
Furthermore (because of the loops), when $\alpha > 5/8$, the two-sided conformal restriction measure is 
not supported on the set of simple curves.
\end {theorem}

This  concludes the proof of the list of results announced at the end of Section \ref {S3}.

Note for instance that for a given value $\kappa$ (in fact $\kappa =6/5$), one constructs the  
two-sided restriction measure with exponent $2$. So ``from outside'', it is the same as the union 
of two Brownian excursions, but they are constructed very differently (in one case, SLE$_{6/5}$ + loops can be viewed as a path, and in the other case, one has the union of two paths).

When $\kappa \in (8/3, 4]$, there is the problem that $\lambda$ becomes negative. In other words, there is no interpretation in terms of cloud of loops, and the martingale $M_t$ is not a priori bounded anymore. In particular, it can not be a conditional probability (conditional probabilities are anyway bounded by one).
The likely scenario is that even though $M_t$ is not bounded, the martingale still converges in $L^1$. In particular, when $\gamma$ itself exits $H$, the derivative term (that goes to zero) beats the exponential term (that goes to infinity). The case $\kappa=4$ is then critical, in the sense, that if one would add a little more loops, then 
the exponential term would win instead. Anyway, the corresponding exponent $\alpha$ is smaller than $5/8$ when 
$\kappa \in (8/3, 4]$ and it does therefore not correspond to a 
(two-sided) conformal restriction measure anymore since it does 
not exist for these values of $\alpha$.

\section {The Brownian loop-soup}

\subsection {Definition}

The Brownian loop-soup in the plane is a Poissonian realization of Brownian loops (with non-prescribed 
length and non-prescribed starting points) in the plane. In fact, it is a random family of {\em unrooted} 
loops (i.e. loops without special marked points on it that can be viewed as its origin). More precisely,
we call an unrooted loop with time-length $T$ a continuous mapping $l$ from 
$T S^1$ onto $\C$, where different loops obtained by shifting time are
identified (i.e. $l(T \exp (i \cdot))$ and $l( T \exp (i (\theta + \cdot))) $ are the same unrooted
loop).

Define the law $P_{x,T}$ on Brownian loops starting and ending at $x$ with time-length $T$. This is simply
Brownian motion started from $x$ and properly conditioned to be at $x$ at time $T$.
It is possible to view this as a probability measure on unrooted loops. 
Then, we use non-deterministic starting points and time. More precisely, we define the infinite 
measure on unrooted loops by
$$
\nu = \int_{\C} dz^2 \int_0^\infty \frac{dT} {T^2} P_{z,T}.
$$
Finally, we define $\nu_D$ to be the measure $\nu$ restricted to those loops that stay inside a 
given domain $D$.

This definition of $\nu_D$ has some very nice properties that recall conformal restriction:
\begin {itemize}
\item
If $D \subset D'$, then $\nu_D$ is equal to $\nu_{D'}$ restricted to those loops that stay in $D$
(this is obvious from the definitions).
\item
If $\Phi$ is a conformal transformation from $D$ onto $D'$, then the image of the 
measure $\nu_D$ under the transformation $\Phi$ is $\nu_{D'}$.
\end {itemize}
This second property is quite strong, but it is easy to prove. Let us stress that in
 order for this property to hold, it is 
important to work with unrooted loops.

\medbreak

The Brownian loopsoup with intensity $\lambda$ is a Poissonian realization of $\lambda \nu$ (or $\lambda \nu_D$): It is a random countable family of unrooted Brownian loops. It inherits the properties of the measure $\nu$:
 
 \begin {itemize}
 \item
 If $(l_j, j \in J')$ is a Brownian loopsoup in $D'$, and if $J = \{ j \in J' \ : \ l_j \subset D \}$ for
 some $D \subset D'$, then $(l_j, j \in J)$ is a Brownian loopsoup in $D'$.
 \item
 If $(l_j, j \in J)$ is a Brownian loopsoup in $D$ and if $\Phi: D \to D'$ is a conformal transformation, then 
 $(\Phi (l_j), j \in J)$ is a Brownian loopsoup in $D'$ (modulo time-changing the loops).
 \end {itemize}
 
\subsection {Relation between loops}

Throughout this paragraph, the ``intensity'' $\lambda$ will be fixed.
Suppose that $\gamma$ is a path from zero to infinity in the upper half-plane (we do not have to assume that it 
is a simple path, it could as well be a path that bounces on its own past i.e. a path with double points but no 
``self-crossings'', or even just a Loewner chain). Suppose that it is parametrized in the way that we described for Loewner chains.  We have already described one way to attach Brownian loops 
to $\gamma$ (the generalization to the intensity $\lambda$ of the way that one would
 attach Brownian loops to SLE$_2$ in order to recover a 
Brownian excursion). This definition is ``dynamic'': For each time $t$, one tries to 
add a Brownian loop in $\H \setminus \gamma [0,t]$ that starts and ends at $\gamma_t$.

Let us now define another way. Suppose that one has a random Brownian loop soup $(l_j, j \in J)$ in the 
upper half-plane (with intensity $\lambda$). For each loop in the 
soup, either it is hit by $\gamma$ or not. For each loop $l_j$ that is hit by $\gamma$, we attach it 
to $\gamma$ (at the first time $t_j$ at which it hits $l_j$).

These two procedures are a priori different, but:

\begin {theorem}
These two procedures to add randomly add loops to $\gamma$ are identical.
\end {theorem}

This has a number of rather unexpected consequences:
\begin {itemize}
\item
In the second construction, the right-boundary of the set obtained by adding the loops to $\gamma$
is identical to the set obtained to adding to $\gamma$ only those loops that do intersect the  
right-boundary of $\gamma$ (this remark is non-trivial only when $\gamma$ is not a simple curve).
Hence, it shows that adding the loops dynamically to $\gamma$, or to its right-boundary
also creates (in law) the same right-boundary. This will be important in the discussion of the 
``duality'' conjectures. 

\item
Adding loops dynamically to $\gamma$ or to the time-reversal of $\gamma$ (i.e. to view $\gamma$ as a path 
from infinity to the origin) is the same (even if a loop does not appear at the same ``time'' in both cases).
This will be important in the discussion of the ``reversibility'' conjectures.

\item
In the previous section, we have seen that when $\kappa \ge 8/3$ adding loops to an SLE$_\kappa$ (and filling) creates a two-sided restriction measure. Hence we have a formula for the probability that there exists a loop 
in the loop-soup that intersects both the SLE and the complement of $H$ (this probability is $\Phi_H'(0)^\alpha$).
This probability is clearly the same as that of the event that the SLE does not hit the set formed 
by the complement of $H$ and all the loops in the loop-soup that intersect the complement of $H$. In other 
words, $\Phi_H'(0)^\alpha$ represents the probability that the SLE does avoid a random set (one attaches the loops
to the complement of $H$ instead of to the SLE).
  
\end {itemize}

\section {SLE($\kappa, \rho$) processes}

\subsection {Bessel processes}

Let us first recall a few elementary facts on a special class of one-dimensional Markov processes: The Bessel processes.
Suppose that we are looking for a non-negative one-dimensional Markov process $X$ started from the origin that has the 
same scaling properties as Brownian motion (i.e. multiplying time by $K^2$ is the same in law as multiplying space by $K$). Under mild conditions, $X$ will be solution of a stochastic differential equation:
$$ dX_t = \sigma (X_t) dB_t + b (X_t) dt .$$
Scaling indicates that $\sigma$ should be constant, and that $b$ should be of the type $cst /x$.
One can scale out the constant $\sigma$ by changing $X$ into $X/ \sigma$, so that one is left
with the SDE
$dX_t =dB_t + c dt/ X_t$.
The solution to this SDE is called the Bessel process with dimension $d= 1+ 2c$. Examples are for instance the 
modulus of $d$-dimensional Brownian motion (that is obviously a Markov process with the right scaling), when $d$ is a
positive integer.
It is easy to see that:
\begin {itemize}
\item
These processes are well-defined and exist if $d>1$. They satisfy for all $t \ge 0$
\begin {equation}
\label {e4}
X_t = B_t + \int_0^t \frac {d-1}{2 X_s} ds 
\end {equation}
(when $d \le 1$, this last equation cannot hold, one has to introduce an additional local time push 
when $X$ hits the origin).
\item
When $d \ge 2$, the process never hits the origin for positive times. It does hit the origin 
infinitely often if $d <2$, but (\ref {e4}) still holds as long as $d>1$.
\end {itemize}

Bessel processes appear in various settings, as soon as a scaling property is combined with a Markov property (for instance for the so-called Ray-Knight Theorems, see e.g. \cite {RY}, where a detailed study of Bessel processes can be found). 

In fact, one can view SLE as a two-dimensional version of the flow generated by the stochastic differential equation
(\ref {e4}). If $g_t$'s are the conformal mappings associated to an SLE$_\kappa$, then define
$$ X_t^z = \frac { g_t (z) - W_t } {\sqrt {\kappa}}.$$
It satisfies
$$
dX_t^z = - dB_t + \frac {2} {\kappa X_t^z} dt.
$$
So, $X_t^z$ can be viewed as  the complex flow of the Bessel process of dimension $1 + (4/ \kappa)$.
The phase transition for SLE at $\kappa =4$ corresponds to the phase transition for Bessel processed at 
$d=2$:
When $\kappa \le 4$, the SLE is almost surely a simple path with no double points, while when $\kappa > 4$
it is a path with many double points (and it hits the real line infinitely often, which corresponds to the 
fact that $X_t^z$ hits almost surely zero for all real $z$).

\subsection {Definition}

One motivation of what follows is to try understand the law of the (simple) right-boundary of 
a restriction measure. When $\alpha =5/8$, we know that is it SLE$_{8/3}$, but what happens when $\alpha > 5/8$? 
By scale-invariance of the restriction measure, the right-boundary should be also scale-invariant.
Recall that the right-boundary should be a simple curve from the origin to infinity. Hence, it can be viewed 
as obtained from a random continuous real-valued driving function $W_t$ via Loewner's equation just 
as SLE$_{8/3}$ is obtained from $\sqrt {8/3}$ times Brownian motion. The question is therefore to understand the 
law of this random driving function $W_t$. The right-boundary $\gamma$ is then defined by 
$$
\partial_t g_t (z) = \frac {2} {g_t (z) - W_t }
$$ 
and $\gamma (t) = g_t^{-1} (W_t)$.

This process $W$ is a priori not Markovian, and does not have stationary increments. But it should have the same 
scaling property as Brownian motion has. Also, it is rather natural to assume that if $O_t$ is the image of the origin 
(more precisely the ``left'' image of the origin) under $g_t$, then  the law of 
$g_t( \gamma [t, \infty))$ depends on the past only via the location of the two points
$O_t$ and $W_t$. Recall that 
$$ 
O_t = \int_0^t \frac {2 ds }{O_s - W_s}
,$$
because $O_t= g_t (0)$ by Loewner's equation. It therefore also satisfies the same Brownian scaling 
as $W$. Hence, we are led to conclude that $W-O$ is Markov and satisfies the Brownian scaling property.
It should therefore be (the multiple of a) Bessel process.

Suppose that $\kappa >0$ and $\rho > -2$ are fixed. We now define $X$ as a Bessel process of dimension 
$1 + 2 (\rho+2) / \kappa$, started from the origin. We will want that the $\sqrt \kappa X = W - O$.
We therefore define
$$ O_t =- 2 \sqrt \kappa  \int_0^t {ds}{ X_s }$$
and
$$
W_t = O_t + \sqrt {\kappa} X_t .
$$
Then,
$$
dW_t = \frac {\rho}{W_t - O_t} dt + \sqrt{\kappa} dB_t.
$$
in other words, the driving process is just as for SLE$_\kappa$, but it gets an 
additional (scale-invariant) push from the left-image of the origin.
This push is repulsive or attractive depending on the sign of $\rho$. 
We call the (two-dimensional) path that is generated by this random driving function an 
SLE($ \kappa, \rho$).
When $\rho=0$, it is ordinary SLE$_\kappa$.
The construction of SLE($\kappa, \rho$) shows that it is scale-invariant.

\subsection {SLE($\kappa,\rho$) and restriction}

It turns out that the previous heuristic guess holds and indeed produces the 
right boundaries of the conformal restriction measures:

\begin {theorem}
The right-boundary of the one-sided conformal restriction measure of exponent $\alpha$ is 
SLE($8/3, \rho$), where 
$$ \rho = \frac { -8 + 2 \sqrt {24 \alpha +1 }} 3$$
(or equivalently, $\alpha = (3 \rho +10) (2 + \rho) / 32$).
\end {theorem}

This is not really surprising; the dimension of the right-boundaries had to be $4/3$ and this 
already forces $\kappa$ to be equal to $8/3$. Then, one is left with the one-dimensional family of 
SLE($8/3, \rho$)'s that are the unique scale-invariant candidates for the one-dimensional family of 
right-boundaries of restriction measures.
In order to prove this Theorem rigorously, one proceeds roughly as in the $\rho = 0$ case:
One just has to find the correct martingales using It\^o's formula. It turns out that the martingale is here
$$
M_t 
= h_t' (W_t)^{5/8}  h_t'(O_t)^b \left(\frac { h_t (W_t) - h_t (O_t) } { W_t - O_t}\right)^c
$$ 
for well-chosen $b$ and $c$
(it is $b=\rho (4 + 3 \rho) /32$ and $c = 3 \rho + 8$), and studying the asymptotic behaviour of this martingale 
shows that it is indeed the conditional probability of 
$ \{\gamma \subset H \}$ given $\gamma [0,t]$, and that $\gamma$ therefore satisfies conformal 
restriction with exponent $ 5/8 + b + c= \alpha$ because
$$
\P [\gamma \subset H ] = E [ M_0] = \Phi_H' (0)^{5/8 + b+c}. 
$$
Recall that a Bessel process of dimension $d$ hits the origin almost surely if and only if $d<2$.
Here, it is not difficult to see that the process $X$ hits the origin if and only if 
the SLE($8/3, \rho$) curve hits the negative half-line.
Hence:

\begin {proposition}
Suppose that the simple curve $\gamma$ satisfies one-sided restriction with exponent $\alpha$.
Then, it does not intersect the negative half-line (almost surely) if and only if $\alpha \ge 1/3$.
\end {proposition}

This was not obvious at all in the ``reflected and conditioned Brownian motion construction'' or in the ``Poissonian
cloud of excursions construction'' of these curves $\gamma$.

\subsection {Interpretation in terms of non-intersection}

The previous martingale can be interpreted  via non-intersection between
independent samples of restriction measures.
We will here just give the hand-waving interpretation. Suppose that one is considering an SLE$_{8/3}$
(from the origin to infinity), and an independent one-sided filled restriction sample $A$ with 
exponent $b >0 $, started also from the origin, and that we condition $\gamma$
and $A$ not to intersect.
Of course, this is an event of probability zero, but it is possible to give to this a rigorous 
meaning by an appropriate asymptotic procedure.

Then, one gets a ``conditional'' joint law for $(A, \gamma)$. The marginal law of $\gamma$ is then 
that of an SLE($8/3, \rho$), where $\rho (4 + 3 \rho) / 32 = b$. In other words, one can interpret the 
repulsive push in the definition of SLE($8/3, \rho$) as a conditioning (not to intersect another
independent restriction measure).
With this interpretation, it is not surprising that the SLE($8/3, \rho$) satisfies one-sided restriction
(the non-intersection property is also ``restriction-invariant''). 

We therefore get the following interpretation of the exponents $b$, $c$ and $\alpha$:

\begin {itemize}
\item
If one considers an SLE$_{8/3}$ and a one-sided restriction measure, the exponent that measures 
how unlikely it is that they do not intersect is $c$ (for instance, the probability that if they 
start at $\eps$ apart, they reach distance one without intersecting is of the order $\eps^c$).
\item
The exponent of the conditioned SLE (not to intersect the one-sided restriction measure of exponent $b$)
is $5/8 + b + c = \alpha$. This is the exponent called $\tilde \xi (5/8, b)$ in the 
papers \cite {LW2,LSW1}.
\item
The law of the conditioned SLE is SLE(8/3, $\rho$).
\end {itemize}

For instance, when $b= 5/8$, then one gets the law of two SLE$_{8/3}$'s conditioned not to intersect. 
We see that $\alpha = 2$. Hence, the law of two SLE$_{8/3}$'s conditioned not to intersect is the same as
the filling of the union of two independent Brownian excursions. Note that the fact that this exponent 
should be $\alpha=2$ can be heuristically understood without using SLE just working with self-avoiding walks (see \cite {LSWSAW}).

One can then iterate this: Condition two $SLE_{8/3}$'s that are already conditioned not intersect 
each other, not to intersect a third one. This leads to another exponent, and to the description of these
conditioned paths. This gives the half-plane intersection exponents that had been predicted by
Duplantier and Saleur \cite {DS}  for self-avoiding walks in a half-space (i.e. surface exponent polymers).
 Note that this approach does not only provide the value of the exponents, but also 
the precise description of the laws of the conditioned paths. 

More generally, we see that all half-space critical exponents that had been derived (via SLE) or predicted 
in theoretical physics (via conformal field theory, Coulomb gas methods or 
the KPZ quantum gravity functions),
seem to fit in the present ``conformal restriction'' framework. In the special case $\kappa=2$, one for instance 
recovers the exponents derived by Rick Kenyon \cite{Kenn} for loop-erased random walks and uniform
spanning trees (these exponents had been predicted in \cite {Dle}) using the link with tiling 
enumerations. 
The other (full-plane) exponents can be similarly worked out using ``radial restriction'' \cite {LSWrr}. 

When $\kappa \not= 8/3$, one can still define SLE($ \kappa, \rho$) in the same way. In this case, the 
martingale $M$ involves an additional term with the exponential of the integral of the Schwarzian 
derivatives. This shows that for $\kappa < 8/3$,  (SLE($\kappa , \rho$) + loops) still 
satisfies (one-sided) restriction, with an exponent that depends on the two parameters $\kappa$
and $\rho$.

It is worth emphasizing that the ``KPZ function'' proposed by Knizhnik, Polyakov and Zamolodchikov \cite {KPZ}
and used by Bertrand Duplantier \cite {Dqg}
 to predict the exact values of these exponents can be interpreted in a 
simple way in the present setup. Basically, the function $U$ (or $U^{-1}$ depending on the 
conventions used) associated to the exponent $b$ in the complex plane is 
the value $\rho$ such that the SLE($\kappa, \rho$) is SLE conditioned to avoid a 
restriction measure of exponent $b$. 
Recall form \cite {LW1} that the existence of the function $U$ can also be derived by simple considerations
based on restriction-type considerations.
 
\subsection {Another construction of the two-sided measure}

We have seen that the right-boundary of a filled set satisfying two-sided restriction 
with exponent $\alpha$ is an SLE($8/3,\rho$) (because it satisfies one-sided restriction).
This raises naturally the question whether it is possible to construct the left boundary, if one conditions on the
knowledge of this right boundary.

Suppose for a moment that $\alpha \ge 2$. 
The previous interpretation of the martingale $M$ shows that the sample of the two-sided restriction measure
with exponent $\alpha$ can be viewed as the filling of the union of two independent 
samples $\gamma$ and $F_0$ of two-sided restriction measures with respective exponents $5/8$ and $\alpha_0$ (which is a well-chosen
function of $\alpha$), that are conditioned not to intersect (and $F_0$ is to the left of $\gamma$).
The conditional law of the SLE$_{8/3}$ becomes SLE($8/3, \rho$) (the law of the right-boundary of $F$). 
But conditionally on this SLE, $F_0$ is simply the random set $F_0$, conditioned to avoid $\gamma$. As $F_0$ 
satisfies restriction, we get (heuristically, but this can be made rigorous, even when $\alpha \in [5/8, 2]$)
the following construction of the two-sided restriction measure:

\begin {itemize}
\item
Define an SLE($8/3, \rho$) curve $\gamma$ with $\rho= \rho (\alpha)$. This will be the right-boundary of $F$.
Let $\Gamma$ denote the connected component of $\H \setminus \gamma$ that has $\R_-$ on its boundary.

\item
Let $\gamma_0$ denote an SLE($8/3, \rho_0$) curve, where $\rho_0=\rho (\alpha_0)$ in $\H$.
Let $\sigma (\gamma_0)$ denote its symmetric image with respect to the imaginary axis (so that 
it satisfies left-restriction). Finally, let $\gamma_-$ denote the conformal image of $\sigma (\gamma_0)$
 under a (conformal) map from $\H$ onto $\Gamma$ that preserves the origin and infinity.
\end {itemize}

\begin {proposition} 
The set $F$ that has $\gamma$ as its right boundary, and $\gamma_-$ as its left boundary, is a 
random filled set that satisfies two-sided restriction with exponent $\alpha$.
\end {proposition}

Recall that a one-sided restriction path hits the real line almost surely if and only if its exponent is 
smaller than $1/3$. Hence, we see that $\gamma_-$ hits $\gamma$ if and only if $\alpha_0 <1/3$.
This shows that:

\begin {corollary}
A two-sided restriction sample has cut points if and only if $\alpha < 35/24$.
\end {corollary}

Again, this was far from obvious with the SLE+loops construction of this two-sided 
restriction measure.

Recall that it is known that the conditioned Brownian motion has cut points, and \cite {BL} that 
the restriction measure with exponent $2$ doesn't (it is the union of two Brownian paths). 
It should  in fact in be possible to determine the Hausdorff dimension of the set of cut points, and show that
it is $2- \xi (\alpha, \alpha)$ in the notation of \cite {LSW2}, when $\alpha \le 35/24$.

Recall also the fact that for all $\kappa< 8$, the SLE curve has (local) cut points, see Beffara \cite {Be0}.

\section {Relation with reversibility and duality conjectures}

\subsection {Reversibility}
It is expected that if one looks at (chordal) SLE from the origin to infinity in 
the upper half-plane, or at an SLE from infinity to the origin in the upper half-plane, they 
trace (in law) the same path modulo time-reversal. 
This should clearly hold if they are the scaling limits of the discrete models from statistical physics 
that they are supposed to correspond to (e.g. \cite {RS} for these conjectures). More precisely: 

\begin {conjecture}
Suppose that $\gamma$ is the path of chordal SLE$_\kappa$ in $\H$ for some fixed $\kappa \le 8$.
Then, the law of (time-reversed, and time-changed) $-1/\gamma$ is identical to the law of $\gamma$.
\end {conjecture}

In some cases (where we know that SLE is the scaling limit of the discrete model), we know that this 
reversibility conjecture holds: $\kappa=6$, $\kappa=2$ and $\kappa =8$.
For $\kappa = 8/3$, we also know that it holds (it is a consequence of Theorem \ref {main1}).
But for all other $\kappa$'s, this is still an open problem. It is known (Oded Schramm, private communication)
that for $\kappa >8$, it does not 
hold.

The relation between SLE and conformal restriction gives a partial result in this direction:
For instance, if $\kappa < 8/3$, if one adds a Brownian loop soup cloud with intensity 
$\lambda$ to $\gamma$ or to $-1 / \gamma$, then one gets a sample of the same conformal restriction 
measure with exponent $\alpha$ (which is unique). Adding the same loop-soup produces the same set (in law), which advocates in favor of the conjecture.

\subsection {Conditioned Bessel processes}

In this subsection, we recall some relevant fact concerning Bessel processes.
It is well-known that the three-dimensional Bessel process can be viewed as a one-dimensional 
Bessel process conditioned to remain positive, see e.g. \cite {RY}. This statement can be formulated 
precisely in different ways:
\begin {itemize}
\item
If one considers a Brownian motion started from $x \in (0,y)$ and conditioned to hit $y$ before $0$,
then it has the same law as a three-dimensional Bessel process started from $x$ and stopped at its 
first hitting time of $y$ (recall that a three-dimensional Bessel process never hits the origin).
\item
If one considers a Brownian motion $W$ started from $x > 0$, and conditioned no to hit the origin 
before time $T$, then the law of $W[0,t]$ (for fixed $t$) converges to that of a three-dimensional
Bessel process started from $x$ (on the time-interval $[0,t]$) when $T \to \infty$.
\end {itemize}
In both these cases, one can then let $x \to 0$ to say that the three-dimensional Bessel process 
(started from the origin) is Brownian motion (started from the origin) and conditioned to remain forever 
positive.
It is also possible (using the Ornstein-Uhlenbeck process $e^{-t/2}W(t)$) to view the scaled Bessel process 
(scaled in the same way)
as the stationary process (corresponding to the first eigenvalue of the corresponding Dirichlet operator) 
which is the Ornstein-Uhlenbeck conditioned to remain forever positive.  

Just in the very same way, one can say that if $d<2$, the $d$-dimensional Bessel 
process conditioned to remain forever positive is the $(4-d)$-dimensional Bessel process.
It is of course not surprising that this conditioned process is a Bessel process, as it is 
a continuous Markov process, and the Brownian scaling is preserved by the conditioning.

Similarly, the interpretation of SLE($\kappa, \rho$) as $SLE_\kappa$ conditioned not to 
not intersect the sample of a one-sided restriction measure
$\alpha$ that we have described should be understood
in a similar fashion. The relation between Bessel processes that corresponds to this 
conditioning is easily worked out via the Girsanov theorem (that shows in general 
how weighting the 
paths in a certain way is equivalent to an explicit change of measure), see \cite {Whid}.  

\subsection {Duality}

Recall that when $\kappa > 4$, the SLE$_\kappa$ curve is not simple anymore. 
In fact, most points of the curve $\gamma[0,t]$ are not on the boundary of the 
unbounded connected component of $\H \setminus \gamma [0,t]$, i.e. $g_t^{-1} (\H)$. 
At time $t$, one can define a right-boundary (resp. left-boundary), as the part of the boundary of 
the unbounded connected component of $\H \setminus \gamma [0,t]$ that is between
$\gamma (t)$ and $\R_+$ (resp. $\gamma(t)$ and $\R_-$). 
The duality conjecture proposed e.g.  by Duplantier \cite {Dqg}  
states that in some sense, this boundary has similar statistical properties as
SLE$_{\kappa'}$, where $\kappa' = 16/ \kappa$.
Of course, it can not be really SLE$_{\kappa'}$ because it is not a path from the origin to infinity, but 
a path between some random point on $\R_+$ (the last one that $\gamma$ did visit before time $t$) 
and the random point $\gamma(t)$. 

There is some evidence for this given by dimension estimates. The probability that a given point $z \in \H$
is in the $\eps$-neighbourhood (for some large and fixed given time) of this right-boundary decays like 
$\eps^{1- 2/ \kappa}$ i.e. like the probability that $z$ is on the $\eps$-neighbourhood
of SLE$_{\kappa'}$ (Beffara, private communication).
 Note that in order to conclude the the dimension of the outer boundary of
SLE$_{\kappa}$ is $1+  2/\kappa$ for $\kappa \ge 4$, one would need second-moment estimates that are not 
proved at this point.

In order to derive an even stronger relation, one would like to find an exact identity in law between
the outer boundary of a process related to SLE$_\kappa$ and a process related to 
SLE$_{\kappa'}$. It turns out that the SLE($\kappa,\rho$) are useful in this respect.
More precisely, the value $\rho= \kappa -4$  is of very special interest: It can be 
proved that (in some appropriate sense), when $\kappa \ge 4$, SLE($\kappa, \kappa-4$) 
is exactly SLE$_\kappa$ conditioned to never hit the negative half line (this is an event of zero 
probability, but it is however possible to make sense of this). This is due to the fact that for some
$d$, the process $W_t - g_t(0)$ is a Bessel process of dimension $d$ for an SLE$_\kappa$.
The conditioning means that this Bessel process never hits zero; it therefore becomes 
a $(4-d)$-dimensional Bessel process; the corresponding curve is the SLE($\kappa, \kappa -4$)
and  (of course), it does not hit the 
negative half-line. Also, it is not difficult to see that it is a transient curve. 

It will be in fact convenient to define SLE$^\sigma$($\kappa,\rho$) as the symmetric image of
an SLE ($\kappa, \rho$) with respect to the imaginary axis. Hence, SLE$^\sigma$ ($\kappa, \kappa -4$) can be viewed
as SLE$_\kappa$ conditioned not to hit the positive half-axis. Define its 
right boundary. It is a random curve from zero to infinity in the 
upper half-plane, and it is clearly scale-invariant (in law).
Dub\'edat \cite {Dub}, based on the restriction properties of the SLE($\kappa, \rho$) processes 
that we will briefly discuss immediately, has
proposed the following:

\begin {conjecture}
The right boundary of an SLE$^\sigma$($\kappa, \kappa-4$) when $\kappa \ge 4$ is 
 an SLE($\kappa', (\kappa'-4)/2)$.
\end {conjecture}

This result is known to hold for $\kappa=6$ (see \cite {LSWr}, we will come back to this in 
next subsection), it is obvious for $\kappa=4$, and for $\kappa=8$, it 
may follow from \cite {LSWlesl}, via Wilson's algorithm \cite {Wi}.
Note that the SLE$^\sigma$ is an SLE ``repelled from the right'' by the positive half-line, while
the second one is an SLE$_{\kappa'}$ that is attracted from the left (because $\kappa' -4 < 0$).

One reason to propose this conjecture is the following result derived in \cite {Dub}:
Let us fix $\kappa \ge 6$, and let $\tilde \gamma$ denote a SLE($\kappa', (\kappa' -4)/2$).
Consider a Brownian loop-soup with intensity $\lambda (\kappa)$. Define the right-boundary of
the set obtained by adding the loops in the loop-soup that it intersects to the SLE$^\sigma$($\kappa,\kappa-4$) $\gamma$. It 
has the same law as the one that one obtains by adding the same loop-soup to $\tilde \gamma$.

Note that even though the reversibility conjecture fails to be true for $\kappa >8$, the duality conjecture should hold even for these values of $\kappa$. It may in fact be the case, that some ``conditioned'' versions of SLE$_\kappa$
when $\kappa > 8$ (for instance ``conditioned not to hit the real line'') do satisfy reversibility.

This leads of course to the following question. The question is interesting on its own, but a 
 positive answer would also solve (at least for $\kappa \le 8/3$) both the 
duality and reversibility conjecture: 

\begin {question}
Suppose that $\gamma^1$ and $\gamma^2$ are two random simple curves from the origin to infinity in the 
upper half-plane. Suppose that $\lambda \ge 0$ is fixed and define the right-boundary $\gamma_+^1$ (respectively $\gamma_+^2$) of the union of $\gamma^1$ (resp. $\gamma^2$) and the loops $l_j$ in a loop-soup of intensity $\lambda$ 
that it does intersect. 
Assume that the law of $\gamma^1_+$ and $\gamma_+^2$ are identical. 
Are the laws of $\gamma^1$ and $\gamma^2$ necessarily identical too?
\end {question}

For related considerations, see \cite {Dub}. Note that reversibility (and duality) would for instance also follow 
if one can prove that these SLE are scaling limits of discrete models for which these facts are satisfied.

\subsection {SLE$_6$ and SLE$_{8/3}$, locality and restriction.}

The random process SLE$_6$ is the unique possible conformally invariant scaling limit 
of critical percolation cluster interfaces (see e.g. \cite {Wstf}), and this has been proved 
by Smirnov \cite {Sm} to hold in the case of critical site percolation on the triangular lattice.
By construction, it is clear that this is not a simple curve. It has (many) double points, and it hits
the real axis a lot. In particular, it can not satisfy restriction.

However, one can (in the same way as before) construct an SLE$_6$ that is ``conditioned not to 
hit the positive real axis''. This is the SLE$^\sigma$(6,2). But, in the discrete case of percolation, it 
is not difficult to see that if one further conditions the interface not to exit a domain $H$, one just has the 
conditioned percolation path conditioned not to intersect the ``positive'' part of the boundary of $H$.
Hence, in the scaling limit, this conditioned SLE$_6$ path should satisfy conformal restriction.
It is therefore not surprising that its right-boundary is an SLE($8/3, \rho$) for some $\rho$.
Note also that the corresponding discrete measure is also an intrinsic measure; the weight of a path 
depends only on its number of neighbours.

Actually, one can also work out another relation between SLE$_6$ and reflected Brownian motion 
with angle $\theta = 2 \pi /3$ (see \cite {W2000,LSWr}), based on their ``locality property''. This 
shows that the (one-sided) restriction exponent of SLE$^\sigma$($6,2$) is $1/3$.
We will not go into this locality property here, but we just mention that locality is related to 
restriction in that if a process satisfies locality, then this process appropriately conditioned not  
to hit the boundary satisfies restriction.

One can in fact also condition the SLE$_6$ to hit neither the left boundary, nor the right boundary.
For the same reasons (due to locality of SLE$_6$), one sees that the obtained path should satisfy 
two-sided restriction. In fact, it is possible to a priori argue that the restriction exponent 
of the obtained path must be one (see the a priori estimates in \cite {LSW5}).

\section {Conformal restriction and representations}

\subsection {Algebraic background}

Define the algebra ${\cal A}$ generated by the vectors $(l_n)_{n \in \Z}$ satisfying the commutation relations
$$ [ l_n, l_m] = (n-m) l_{n+m}.$$
Because of these relations, it is easy to see that a basis of this algebra is given by the family of 
vectors of the type $l_{n_1} l_{n_2} \ldots l_{n_p}$ for
$n_1 \le n_2 \le \cdots \le n_p$. This algebra is often called the algebra of polynomial vector fields
on the unit circle, because it can be realized as $l_n = - z^{n+1} d/dz$.

Suppose that there exists a vector space $V$ on which the algebra ${\cal A}$ acts i.e., one view ${\cal A}$ as
a subalgebra of the set of endomorphisms of $V$. If it happens that for some vector 
$v \in V$, 
$ l_n v = 0 $ for all $n >0$ and $l_0 v = hv$ for some real $h$, then we say that this is a highest-weight 
representation of ${\cal A}$. It is then possible to say that 
${\cal A}$ acts only on the vector space generated by 
$v$ and its ``descendants''  of the type $l_{n_1} \ldots l_{n_p} v$, where 
$n_1 \le \cdots \le n_p < 0$ (because $l_n v =0$ as soon as $n >0$).

In fact, such representations can be constructed for all values of the highest weight $h$. 
The vector space is ``naturally'' graded i.e. it can be decomposed as the direct sum of 
$V_n$'s for non-positive $n$'s, where $V_n$ consists of all vectors $w$ in $V$ such that 
$l_{-n} w$ is co-linear to $v$. For instance, $V_{-1}$ is generated by $l_{-1}v$, 
$V_2$ is generated by $l_{-2}v$ and $l_{-1}^2 v$, and more generally, $V_{n}$ is generated by the 
family $F_n$ of vectors $l_{n_1} \ldots l_{n_p} v$, where $n_1  \le \cdots \le n_p < 0$ and $n_1 + \cdots +
n_p = n$. Then, each $l_m$ maps $V_n$ onto $V_{n+m}$.

The representation is said to be degenerate at level $n$, if in fact the vectors of $F_{-n}$
are not independent. For instance, it is degenerate at level 2 if 
$l_{-1}^2 v$ and $l_{-2}v$ are in fact co-linear.
Suppose for instance that for some $\kappa$,
$$\kappa l_{-1}^2 v = 4 l_{-2} v.$$
Then, applying $l_2$ to both sides, and applying the commutation relations, we see that 
$$ \kappa \times 3 \times 2  \times l_0 v    = 4 \times 4 \times   l_0 v$$ 
Hence, if we assume that $h \not= 0$, we see that $\kappa$ has to be equal to 
$8/3$.
If we apply $l_1$ to both sides instead of $l_2$, we get that
$$ \kappa (2l_0 l_{-1} + l_{-1} l_1 l_{-1}) v = \kappa (2 l_{-1} h v + 2 l_{-1} v + l_{-1} 2 h v ) =
  12 l_{-1} v$$ 
so that (unless $l_{-1}v =0$), 
$h=5/8$.
It is in fact indeed possible to construct such a representation that is degenerate at level two
(with $\kappa = 8/3$ and $h=5/8$).

\subsection {Relation with restriction}

Suppose that $\gamma$ is a one-sided restriction curve with exponent $\alpha$.
It is in fact possible to construct a highest-weight representation of ${\cal A}$ associated
to $\gamma$, and the corresponding highest-weight $h$ is just the exponent $\alpha$.
One proceeds as follows:
Define for all positive $x_1, \ldots, x_N$ the renormalized probability the the 
curve $\gamma$ passes in the neighbourhood of these $N$ points. 
More precisely,
$$
B_N (x_1, \ldots, x_N) = \lim_{\eps \to 0} \eps^{-2N}
P [\gamma \cap [x_j, x_j + i \eps ] \not= \emptyset, j =1 , \ldots, N] 
.$$
The fact that $\gamma$ satisfies one-sided restriction implies a certain relation between 
$B_N$ and $B_{N+1}$: Suppose that one considers an $N+1$-th infinitesimal slit $[x,x+i \eps]$.
Then, either the path $\gamma$ hits it also (and this probability is given by $B_{N+1}$) 
or it avoids it (and the probability that it hits the $N$ other ones is now given in terms of
$B_N$ and the conformal mapping from $\H \setminus [x ,x +i \eps]$  onto $\H$.
This relation can be written as:
$$
B_{N+1} (x,x_1, \ldots , x_N)
= \frac{\alpha}{x^2} B_N (x_1, \ldots, x_N)
+ \sum_{n\ge 1} x^{n-2} L_{-n} B_{N} (x_1, \ldots, x_N)
$$
for some operators $L_N$. This equation can then be rephrased in terms of a highest-weight 
representation of ${\cal A}$ with highest-weight $\alpha$, see \cite {FW} for details. Basically,
one shows that these $L_n$'s, when defined on appropriate functions, do satisfy the same commutation
relation as the $l_n$. 

If one supposes that the curve $\gamma$ satisfies also $({\cal P})$. Then, this leads 
to some martingales that describe the conditional (renormalized) probability to hit the
infinitesimal slits. After a very small time $t$, the new conditional probability is 
roughly
$$
g_t'(x_1)^2 \ldots g_t'(x_N)^2 B_N ( g_t (x_1) - W_t, \ldots, g_t (x_N) - W_t )
$$
(the derivative terms are due to the fact that the sizes of the slits vary).
Hence, by It\^o's formula, one sees readily that for the above-mentioned representation
$$\kappa / 2 l_{-1}^2 v - 2 l_{-2} v = 0
$$
i.e. that it is degenerate at level two. This explains why the same values $\kappa =8/3$ and
$h=   \alpha = 5/8$ show up.

A more involved study can be applied to the case where the right-boundary $\gamma$ is constructed 
via an SLE to which one adds loops, as described before. Again, it is possible to recover the 
relation between the density of loops $\lambda$, the parameter $\kappa$ and the highest-weight (or exponent)
$\alpha$ from algebraic considerations. This relation is the same as the one, that one obtains when studying 
highest-weight representation of the Virasoro algebra (the central extension of $V$) that are degenerate 
at level $2$. The quantity $-\lambda$ is then interpreted as the central charge of this representation.
See \cite {FW} for more details. 

Note that in this setup, one ends up with highest-weight representations of ${\cal A}$ itself rather than with 
degenerated highest-weight representations (at level two) of the 
Virasoro algebra, but there is a simple correspondence between them.
  
The fact that the Virasoro Algebra's degenerate highest-weight representations are related to two-dimensional systems lies at the roots of conformal field theory, see e.g. \cite {BPZ0,BPZ,Ca0,Cabook,CaChuo}, and has been one of the ideas that led to the prediction of the exact values of the relevant critical exponents in last decades.
In the recent series of papers \cite {BB1,BB2,BB3,BB4}, Michel Bauer and Denis Bernard have been studying various aspects of the interplay between the conformal field theory and SLE. In particular, they exhibited relations between SLE and such representations in this setup. In their approach, it also turns out \cite {BB4} that as in \cite {FW1,FW, FK}, an instrumental role is played by the local martingale $M$ defined in (\ref {themartingale}).
 
\section {Remarks}

\noindent
{\bf What other discrete measures?}
One can ask the question whether other intrinsic simple discrete measures on paths, will be conformally invariant 
in their scaling limit. A first example goes as follows: Consider the law of simple random walk from
$A$ to $B$ in $D$ (the same as the one used to define $P_{D,A,B}^{BM}$ for bounded $D$) but conditioned 
to have no cut point. It is reasonable to believe that the limiting measure in the scaling limit will in some sense be 
$P_{D,A,B}^{BM}$ conditioned to have no cut-point, which should therefore be conformally invariant (``having no cut-points'' is a conformally invariant property). It is possible (at least on heuristic level) to see
\cite {Whid} that the 
restriction exponent of the limiting measure will be $2$ (this is also related to B\'alint Vir\'ag's \cite {V} 
Brownian beads exponent, and to Beffara's cut time exponents \cite {Be0}).

What happens if one considers the measure on discrete random walks, but this time conditioned to have no triple points.
Will it degenerate? Possibly, one will have to take another penalization (than $4^{-n}$) in order to have a non-trivial
limit. It is likely to be related to a two-sided restriction measure.
Similarly, what happens if one allows no point of multiplicity $k$ instead for $k>3$?

If one penalizes the energy according to the number of double points, does one indeed destroy conformal invariance
as one might at first sight think? Note that the (conditioned) critical percolation interface can be viewed as a measure on paths 
with double but no triple points with a well-chosen intrinsic way of weighting paths.

\medbreak
\noindent
{\bf All is Brownian.}
In some sense, the conformal restriction approach shows that it is probably possible to characterize completely 
SLE (at least for $\kappa \le 8/3$) in terms of planar Brownian motion. The restriction measures can be constructed 
using Brownian motions (conditioned (reflected or not) Brownian motions, for instance). Adding Brownian loops to 
a path is in some sense the unique conformally invariant way to enlarge a given path. And SLEs are (probably) the unique measures on simple paths, such that if one adds Brownian loops, one gets a (Brownian) restriction measure.

While the definition of SLE via iterations of independent identically distributed conformal mappings 
is difficult to generalize to define interfaces on Riemann surfaces, this ``Brownian'' approach seems
well-suited (recall that it is no problem to define Brownian motions on surfaces).
See \cite {Dub2,FK,Z} for progress in this direction.

\end {document}